\documentclass[10pt]{article}

\usepackage{amssymb, amsmath, amsthm, mathtools, color}
\usepackage{multirow}
\usepackage{rotating}
\usepackage{subfigure}
\usepackage{epsfig}
\usepackage{url}
\usepackage{soul}
\usepackage{cite}
\usepackage{graphicx}
\usepackage{mathdots}
\usepackage[english]{babel}
\usepackage{csquotes}
\usepackage{mathrsfs}
\usepackage{overpic}
\usepackage{appendix}
\usepackage{algorithmic} 
\usepackage{algorithm}

\soulregister\cite7
\soulregister\ref7
\soulregister\pageref7
\soulregister\eqref7

\newtheorem*{remark}{Remark}

\newcommand{\bA}{\mathbf{A}}
\newcommand{\bS}{\mathbf{S}}

\newcommand{\bI}{\mathbf{I}}

\newcommand{\overbar}[1]{\mkern 1.5mu\overline{\mkern-1.5mu#1\mkern-1.5mu}\mkern 1.5mu}

\begin{document}
	\title{A Fast and Accurate Algorithm for Spherical Harmonic Analysis on HEALPix Grids with Applications to the Cosmic Microwave Background Radiation}
	\author{Kathryn P. Drake\footnote{Department of Mathematics, Boise State University, 1910 University Drive, Boise, ID 83725-1555 (kathryndrake@u.boisestate.edu, gradywright@boisestate.edu).} $^{,}$\footnote{Corresponding author.} \, and Grady B. Wright$^*$}
	\date{\today}
	
	\maketitle
	
	\begin{abstract}

The Hierarchical Equal Area isoLatitude Pixelation (HEALPix) scheme is used extensively in astrophysics for data collection and analysis on the sphere. The scheme was originally designed for studying the Cosmic Microwave Background (CMB) radiation, which represents the first light to travel during the early stages of the universe's development and gives the strongest evidence for the Big Bang theory to date.  Refined analysis of the CMB angular power spectrum can lead to revolutionary developments in understanding the nature of dark matter and dark energy. In this paper, we present a new method for performing spherical harmonic analysis for HEALPix data, which is a central component to computing and analyzing the angular power spectrum of the massive CMB data sets. The method uses a novel combination of a non-uniform fast Fourier transform, the double Fourier sphere method, and Slevinsky's fast spherical harmonic transform~\cite{FastSPH}.  For a HEALPix grid with $N$ pixels (points), the computational complexity of the method is $\mathcal{O}(N\log^2 N)$, with an initial set-up cost of $\mathcal{O}(N^{3/2}\log N)$.  This compares favorably with $\mathcal{O}(N^{3/2})$ runtime complexity of the current methods available in the HEALPix software when multiple maps need to be analyzed at the same time. Using numerical experiments, we demonstrate that the new method also appears to provide better accuracy over the entire angular power spectrum of synthetic data when compared to the current methods, with a convergence rate at least two times higher.

\end{abstract}

\section{Introduction}
\label{intro}
About 379,000 years after the universe began, the dense plasma of matter cooled enough for neutral hydrogen to form. During this epoch of recombination, the universe was becoming increasingly transparent to photons, which eventually began to move freely through space. Now faintly glowing at the edge of the observable universe, these photons form the Cosmic Microwave Background (CMB) radiation, which has become the strongest evidence for the Big Bang Theory to date~\cite{WMAP1}. While the CMB has been deemed ``the most perfect black body ever measured in nature"~\cite{CMBAnisotropies}, there are temperature and polarization fluctuations that give insight into the primordial universe~\cite{CMBDat}. These anisotropies are consequences of the initial density distribution of matter, and analyzing them can provide a better understanding of the geometry and composition of the universe~\cite{WMAP1,WMAP2}.
\begin{figure}[ht]
	\begin{minipage}{0.4\textwidth}
		\centering
		\includegraphics[scale=.38]{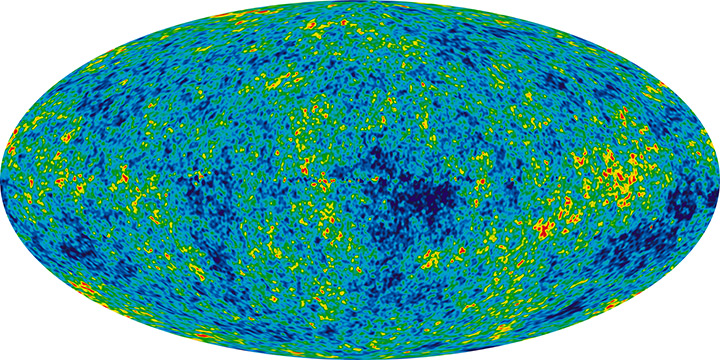}\\
		 \qquad \qquad \qquad \qquad(a)
	\end{minipage}
	\begin{minipage}{0.8\textwidth}
		\centering
		\includegraphics[scale=.4]{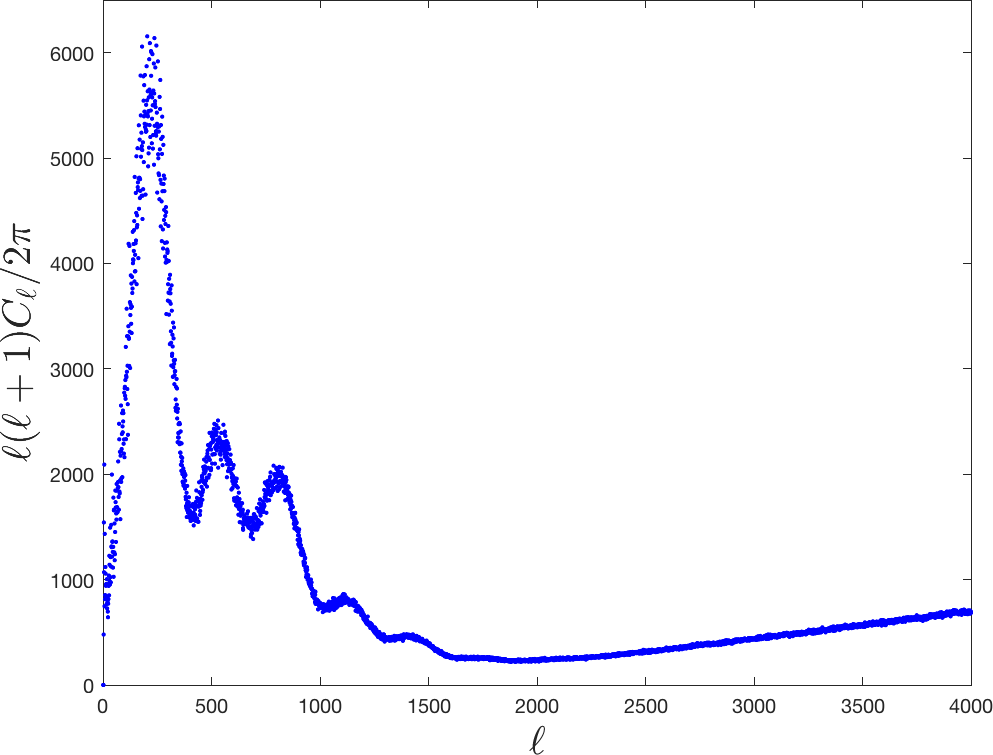}\\
		\quad \,(b)
	\end{minipage}
	\caption{CMB component map from the Planck mission~\cite{PLANCK} (a) and corresponding (scaled) angular power spectrum (b).}
	\label{CMB}
\end{figure}

Using ground-based telescopes, balloons, and satellites which probe the sky in the microwave and infra-red frequencies, scientists have measured the CMB temperature differences at small angular scales. These measurements are quantized and stored as a high resolution sky map of the CMB using the Hierarchical Equal Area isoLatitude Pixelation (HEALPix) scheme~\cite{HEALPA} for the sphere; see Figure~\ref{CMB}a) for an example sky map. Once a sky map is composed, it can then be analyzed by its angular power spectrum. This quantity measures the amplitude of the CMB temperature fluctuations as a function of angular scale and is used to estimate parameters of the cosmological model for the universe~\cite{CMBAnisotropies}. For example, the confirmation of the first peak in the temperature angular power spectrum affirmed that the universe is spatially flat~\cite{CMBHu}. The values of the temperature angular spectrum at higher frequencies are crucial to many aspects of modern cosmology, including the density of dark matter and dark energy in the universe.  The CMB power spectrum (Figure~\ref{CMB}b) is calculated from the spherical harmonic coefficients, $a_\ell^m$, of the sky map as follows:
\begin{equation}
C_\ell = \frac{1}{2\ell+1}\sum_m |a_\ell^m|^2.
\label{APS}
\end{equation}
The spherical harmonic conventions used in this work are detailed in Appendix~\ref{app:sph}. 

The HEALPix scheme~\cite{HEALPA} and the associated eponymous software~\cite{HEALPP} have a number of desirable properties for data collection on the sphere.  First, each pixel has the same surface area and the pixel centers (points) are quasi-uniformly distributed over the sphere.  This is important since any white noise produced by the microwave receivers is exactly integrated into white noise in the pixel area.  Second, the pixels produced by the scheme are based on a hierarchical subdivision of the sphere, which allows for data locality in computer memory and fast search procedures.  Finally, the pixel centers are isolatitudinal, allowing for a significant reduction in the computational cost of performing discrete spherical harmonic transforms---a central component to computing and analyzing the angular power spectrum of the CMB data sets, which from the Planck mission consist of millions of pixels~\cite{PLANCK}.  These properties have made the HEALPix scheme popular for other applications in astrophysics/astronomy~\cite{PhysRevD.93.024013,Majeau_2012,galaxies}, and to several other disciplines, including geophysics~\cite{geodesy}, planetary science~\cite{lunar}, nuclear engineering~\cite{nuclear}, and computer vision~\cite{computervision}.

In this paper, we focus on an aspect of the HEALPix scheme that has received very little attention in the literature: accuracy and computational complexity improvements of the discrete spherical harmonic transform.  We first review the current techniques used in the HEALPix software~\cite{HEALPP}, which are based on equal-weight quadrature, ring-weight quadrature, and pixel-weight quadrature.  We then introduce a new algorithm for computing spherical harmonic coefficients for data collected on HEALPix grids. The main motivation for the method is Slevinsky's recently developed fast spherical harmonic transform (FSHT)~\cite{FastSPH}, which converts bivariate Fourier coefficients for data on the sphere to spherical harmonic coefficients of the data with near optimal complexity. By combining the nonuniform fast Fourier transform (NUFFT)~\cite{NUFFT} and the double Fourier sphere (DFS)~\cite{1sphere} methods, we give a fast and accurate method for obtaining the bivariate Fourier coefficients for functions sampled on the HEALPix grid, which we then use with the FSHT to obtain the spherical harmonic coefficients.  For a HEALPix grid with $N$ pixels (points), the computational complexity of the method is $\mathcal{O}(N\log^2 N)$, with an initial set-up cost of $\mathcal{O}(N^{3/2}\log N)$, which compares favorably with the complexity of the current methods available in the HEALPix software when multiple maps need to be analyzed at the same time. Using numerical experiments, we demonstrate that the new method also appears to be more accurate than the current methods for synthetic data over the whole spectrum, with a convergence rate at least two times higher.  We believe this new scheme will be useful not only for CMB analysis, but also for the many applications of the HEALPix scheme given above that require a spherical harmonic analysis.  Additionally, the algorithm presented here has natural generalizations for other ``equal-area'' isolatitudinal sampling strategies for sphere that do not have a natural way to do fast and accurate spherical harmonic transforms~\cite{SK97,CT98,Leo06,Malkin19}.

The remainder of the paper is structured in the following manner. In section~\ref{HEAL}, we offer supporting information on the HEALPix grid as well as details and analysis of the current methods used in the HEALPix software for computing the spherical harmonic coefficients of CMB maps.  We present the new algorithm for fast spherical harmonic analysis of data collected on the HEALPix grid in section~\ref{HP2SPH}. Numerical results comparing the presented method with that of the methods in the HEALPix software for calculating the angular power spectrum of functions on the sphere are given in section~\ref{NumRes}. Finally, in section~\ref{Conc}, we give some brief conclusions and remarks on future directions of the work.\\

\section{Background and Current Approach}
\label{HEAL}

\subsection{HEALPix Scheme}
\label{HEALP}
The HEALPix scheme\footnote{The HEALPix scheme produces a grid consisting of a collection of pixels of different shapes but the same area.  However, for our method we do not exploit this fact and simply treat the center of each pixel as a point with the given value of the pixel.} was created to discretize functions on the sphere at high resolutions. In addition to creating an equal area pixelization of the sphere, one of the primary motivations behind the scheme was to allow for more computationally efficient spherical harmonic transforms on increasingly large CMB datasets~\cite{HEALPA}. While there are many options for discretizing the sphere, there is no known deterministic method that gives a set of quasiuniform points and allows for exact spherical harmonic decompositions of band-limited functions using equal-weight quadrature. While the HEALPix scheme does not offer optimal complexity for spherical harmonic analyses, it does achieve some efficiency gains over existing schemes for discretizing the sphere. This improvement is accomplished primarily by the isolatitudinal distribution of pixels.

\begin{figure}[h!]
\centering
\begin{tabular}{cccc}
\includegraphics[width=0.23\textwidth]{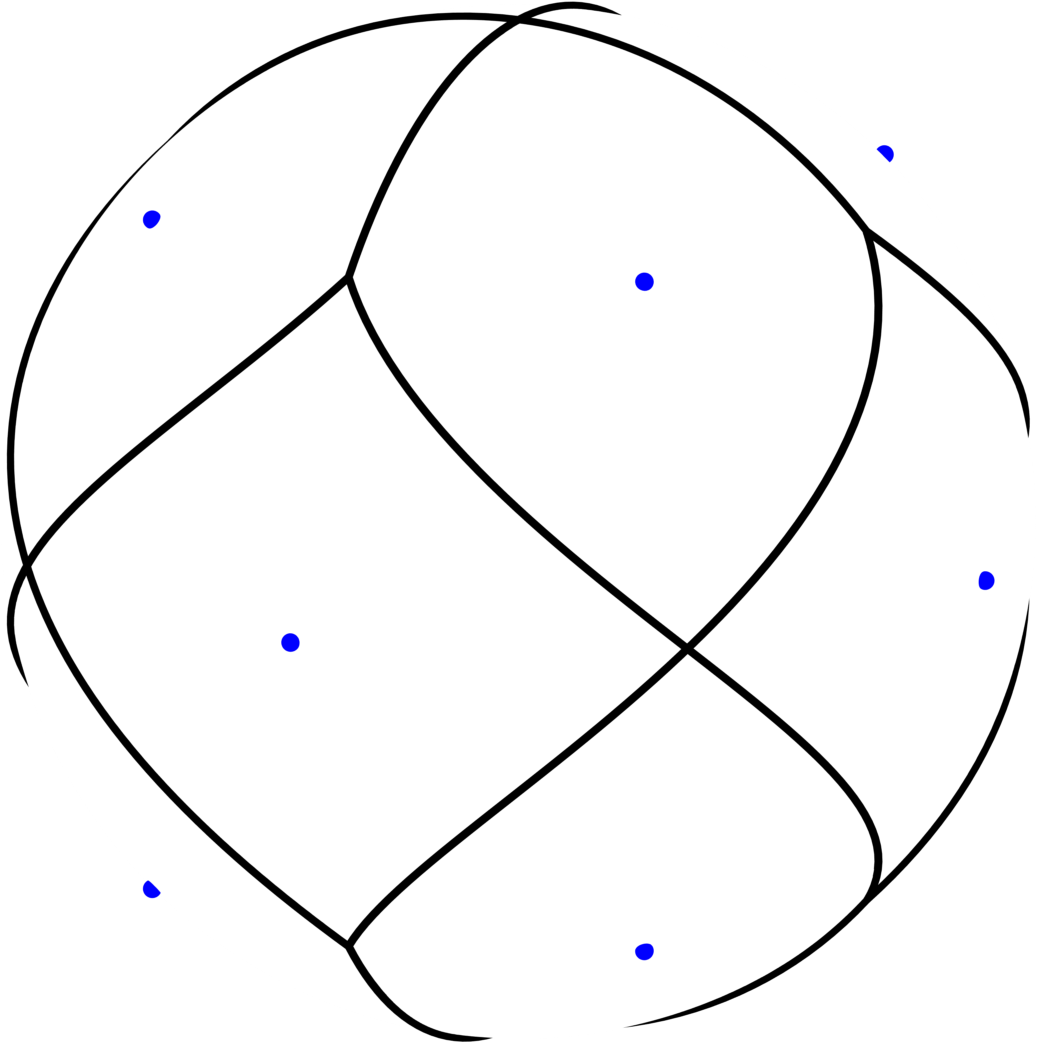} & 
\includegraphics[width=0.23\textwidth]{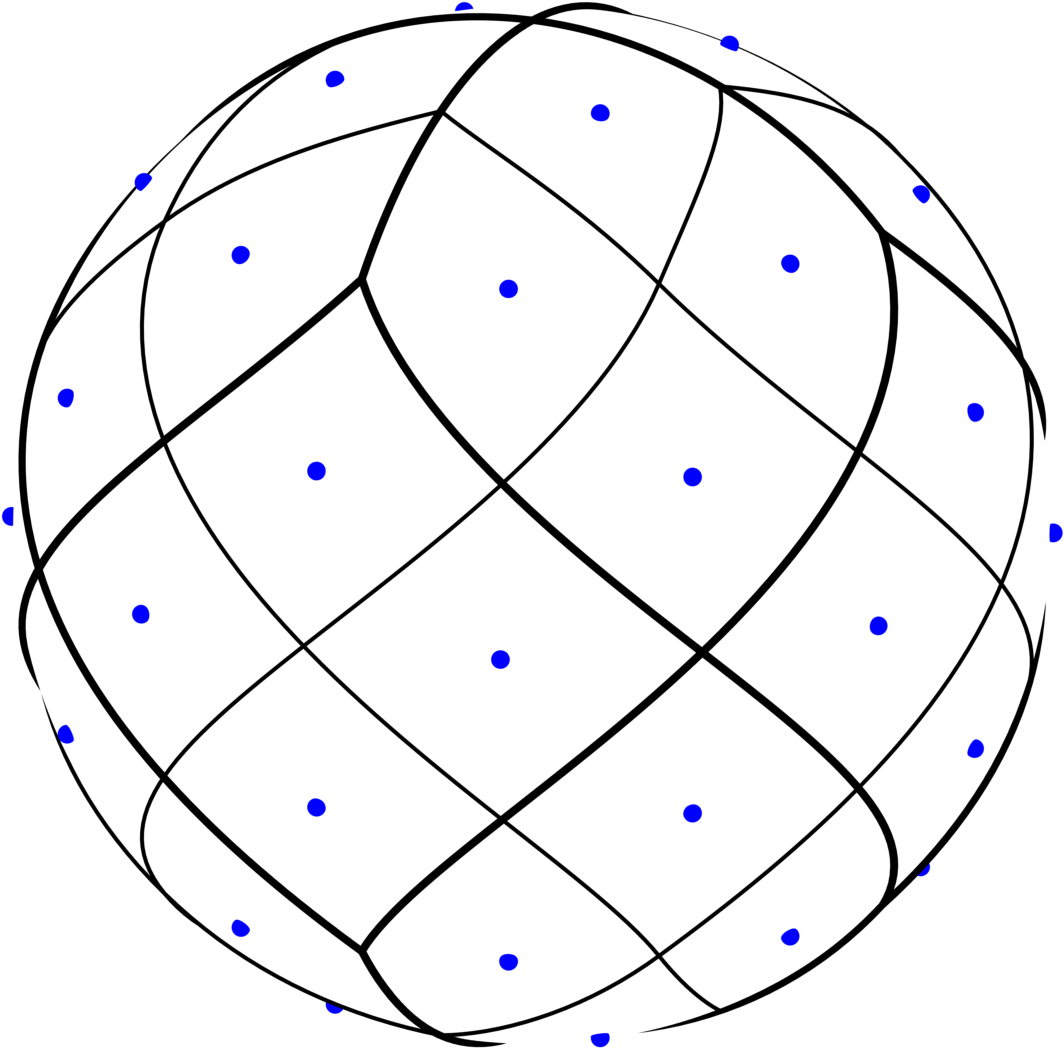} & 
\includegraphics[width=0.23\textwidth]{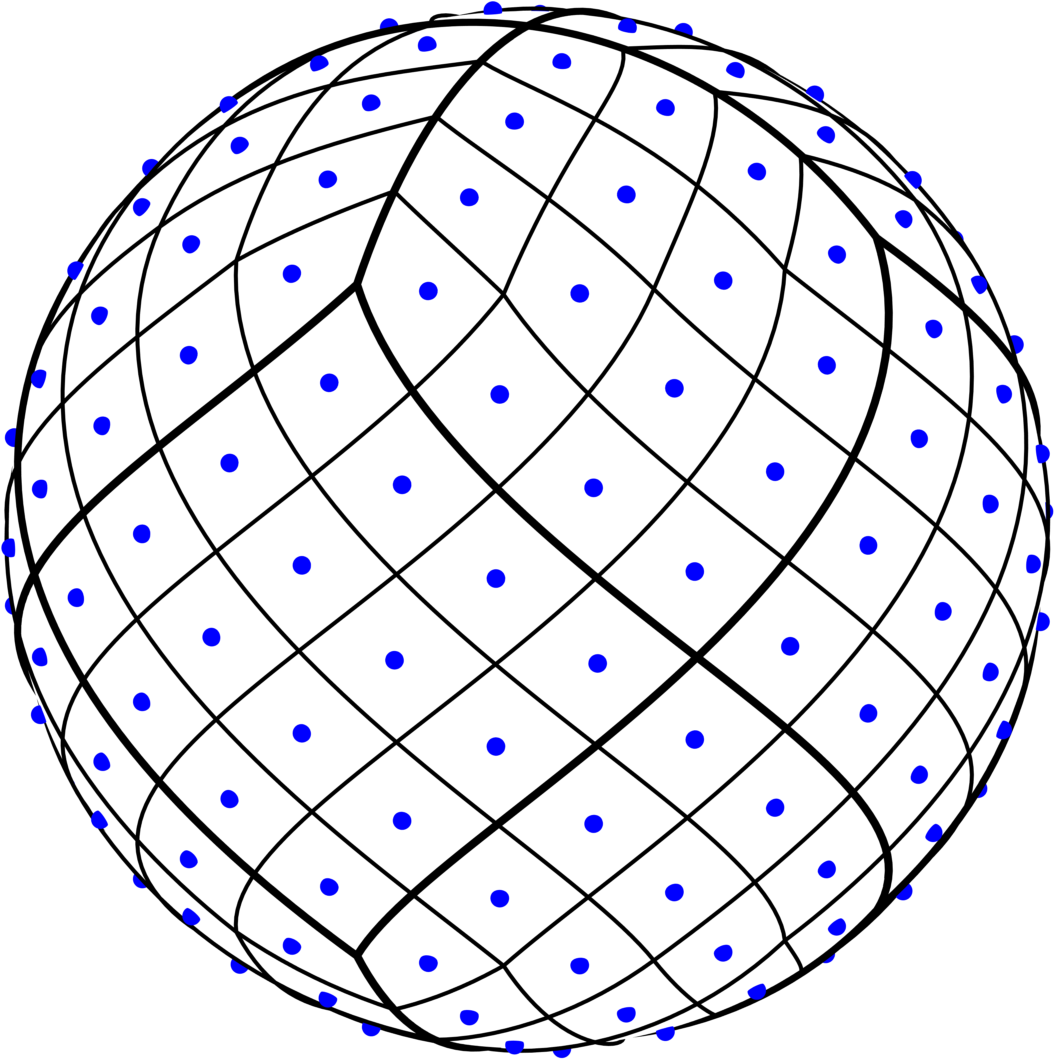} & 
\includegraphics[width=0.23\textwidth]{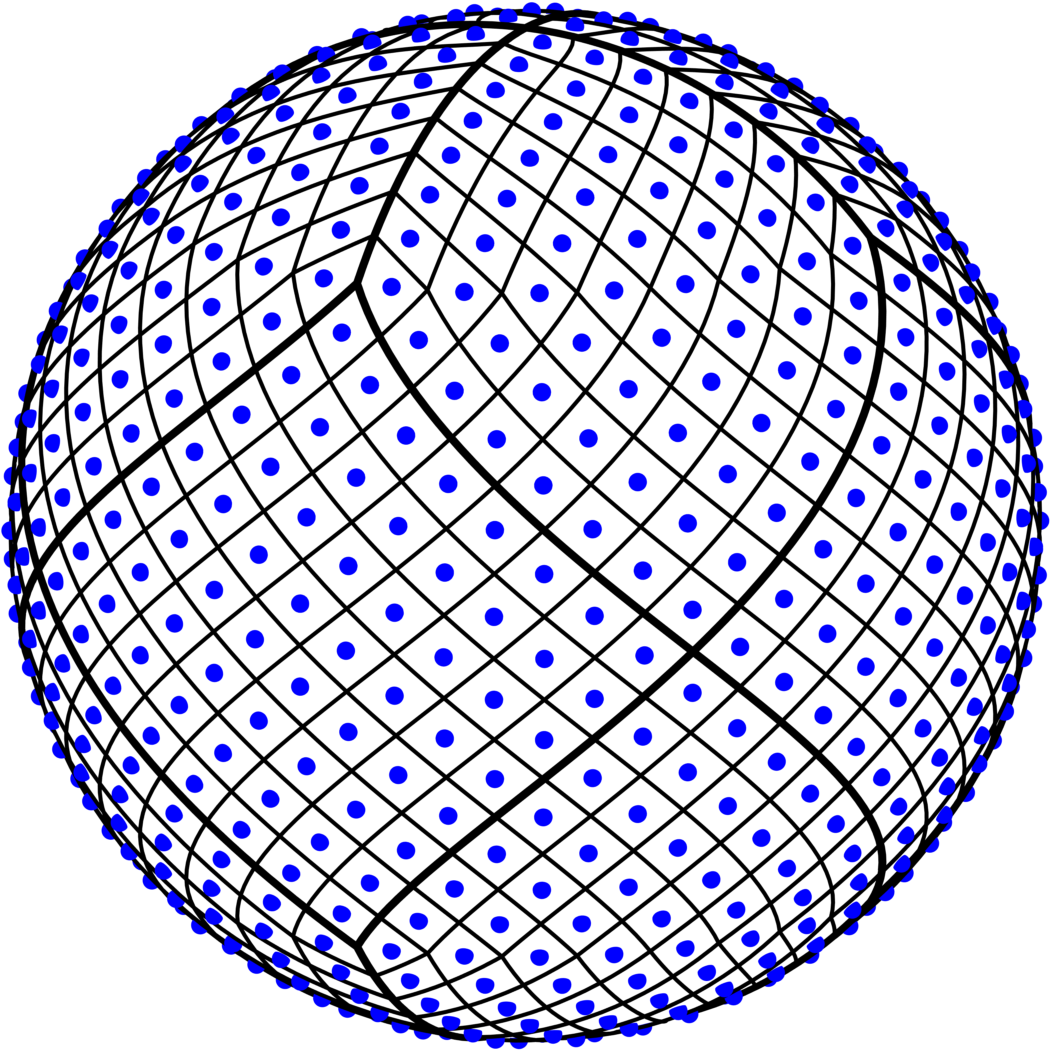} \\
\end{tabular}
\caption{HEALPix grid with resolutions, from left to right, $N_{side}=1,2,4,8$.  The lines indicate the pixel boundaries and the solid dots represent the pixel centers or points.}
\label{Heal_Base}
\end{figure}

\begin{figure}[ht]
\centering
\includegraphics[scale=.45]{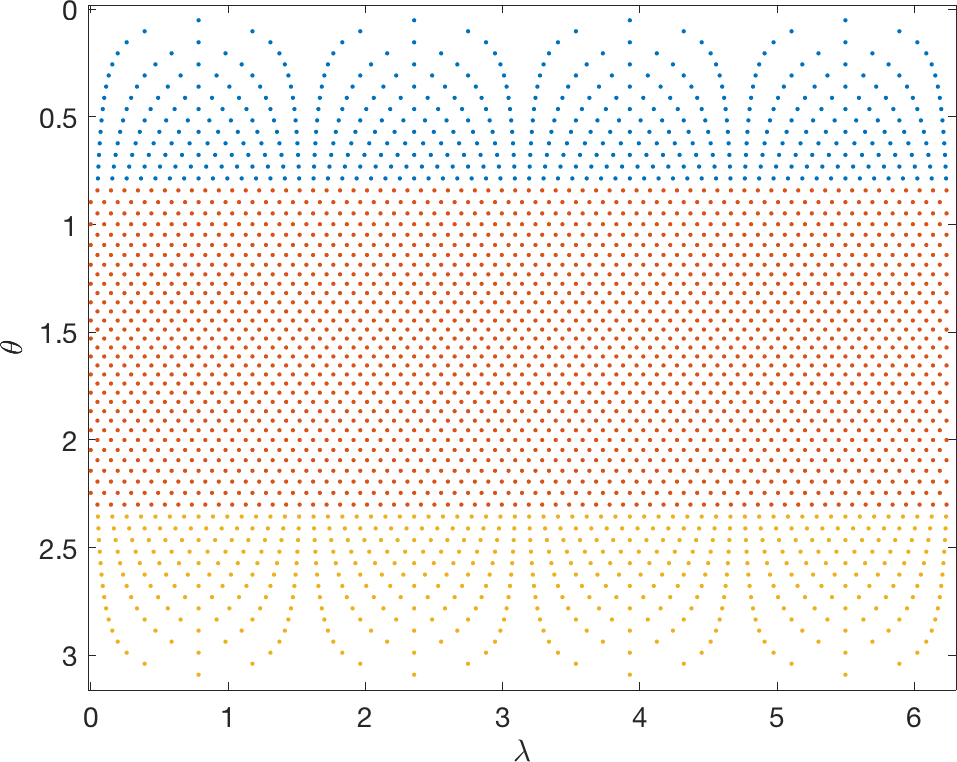}
\caption{HEALPix grid on $[0,2\pi]\times [0,\pi]$, where $\theta$ is latitude, and $\lambda$ is longitude. The point sets in the northern ($\mathcal{N}$), equatorial ($\mathcal{E}$), and southern ($\mathcal{S}$) regions are shown in blue, red, and yellow, respectively.}
\label{Heal_Grid}
\end{figure}

The HEALPix grid resolution is defined using the parameter $N_{side}=2^t, \, t\in\mathbb{N}$, which creates $N_{side}^2$ equal area divisions of each base pixel. Figure~\ref{Heal_Base} illustrates the base resolution grid, $t=0$, and the increasing levels of refinement $t=1,2,3$, where each base pixel is subdivided further into four equal area pixels. A HEALPix map therefore has $N=12 N_{side}^2$ equal area (but differently shaped) pixels, with the centers placed on $4N_{side}-1$ rings of constant latitude. For any $N_{side}$, the HEALPix centers, which we henceforth call the HEALPix points, are defined algebraically using three regions of the sphere, two polar ($\mathcal{N}$ and $\mathcal{S}$) and one equatorial ($\mathcal{E}$)~\cite{Potts}.  In spherical coordinates, the points in these regions are given as
\begin{gather}
\begin{aligned}
\mathcal{N} &:= \left\{\left(\arccos\left(1-\frac{j^2}{3N_{side}^2}\right),\frac{\pi\left(k+\frac{1}{2}\right)}{2j}\right) \; : \; j=1,\dots,N_{side}-1,\; k=0,\dots,4j-1\right\}\\
\mathcal{E} &:= \left\{\left(\arccos\left(\frac{2(2N_{side}-j)}{3N_{side}}\right),\frac{\pi\left(k+\frac{(j+1)\,\textnormal{mod}\, 2}{2}\right)}{2N_{side}}\right) \; : \; j=N_{side},\dots,3N_{side}, \; k=0,\dots,4N_{side}-1\right\}\\
\mathcal{S} &:= \left\{\left(\arccos\left(-\left(1-\frac{j^2}{3N_{side}^2}\right)\right),\frac{\pi\left(k+\frac{1}{2}\right)}{2j}\right) \; : \; j=1,\dots,N_{side}-1,\; k=0,\dots,4j-1\right\}.
\end{aligned}
\label{healpts}
\end{gather}
The final HEALPix point set is $\mathcal{X} = \mathcal{N}\bigcup\nolimits\mathcal{E}\bigcup\nolimits\mathcal{S}$.
The number of points on each ring varies in the polar regions, with only four points on the rings closest to the north and south poles of the sphere, whereas the rings in the equatorial region have the same number of points. This point distribution is illustrated more clearly in Figure~\ref{Heal_Grid}, where the HEALPix points are displayed to a latitude-longitude grid. 

The biggest computational advantage for spherical harmonic analysis in the HEALPix scheme lies in the equally-spaced points on each ring of constant latitude. While this aides computation in the longitude direction with FFTs, the misaligned and unequally spaced points in latitude make fast bivariate Fourier analysis impossible without modification. We address this in the new algorithm presented in section~\ref{HP2SPH}.

\subsection{HEALPix Software Spherical Harmonic Analysis}
\label{HEALsoft}
The standard method in the HEALPix software~\cite{HEALPP} for estimating the angular power spectrum~\eqref{APS} of data at the HEALPix points approximates the exact spherical harmonic coefficients ($\widetilde{a}_\ell^m$) of the data as\\
\begin{equation}
a_\ell^m=\frac{4\pi}{N}\sum_{i=1}^N\overbar{Y}_\ell^m(\lambda_i,\theta_i)f(\lambda_i,\theta_i), \; 0 \leq \ell \leq \ell_{max}, -\ell \leq m \leq \ell,
\label{healsph}
\end{equation}
where $(\lambda_i,\theta_i)$ are HEALPix points in longitude-latitude, $f$ is the data, and $Y_{\ell}^m$ is a spherical harmonic of degree $\ell$ and order $m$ (see Appendix A for a discussion of the spherical harmonic conventions used in this paper).
While the user can input any band limit $\ell_{max}$ for this approximation, the software default is   $\ell_{max}=3N_{side}-1$. Due to the isolatitudinal nature of the HEALPix points, this computation is done with $\mathcal{O}(N^{3/2})$ complexity as opposed to $\mathcal{O}(N^{2})$~\cite{HEALPA}. Note that $N=\mathcal{O}(\ell_{max}^2)$, so the complexity of the $a_\ell^m$ computation is equivalent to $\mathcal{O}(\ell_{max}^3)$. The expression~\eqref{healsph} is a low-order approximation to the continuous inner product~\eqref{sphip} which defines the coefficients, since it uses a simple equal weight quadrature. To improve this approximation, the software employs an iterative procedure, which is referred to as a ``Jacobi iteration''~\cite{HEALPA}.  In order to illustrate how the iterative method converges, we explain it below in the language of linear algebra.

The analysis operation, defined in~\eqref{healsph}, produces an approximation to the spherical harmonic coefficients from the data $f$ on the sphere, whereas the synthesis operation reconstructs the data given the spherical harmonic coefficients:
\begin{equation}
\widehat{f}(\lambda_i,\theta_i) = \sum_{\ell=0}^{\ell_{max}} \sum_{m=-\ell}^\ell a_\ell^mY_\ell^m (\lambda_i,\theta_i),\; i=1,\ldots,N
\label{synthesis}
\end{equation}
Note that we use a hat on $f$ to indicate that computing the spherical harmonic coefficients according to \eqref{healsph} and using them in \eqref{synthesis} gives different function values in general. In matrix-vector notation, we denote \eqref{healsph} and \eqref{synthesis} as
\begin{align*}
\text{Analysis:}\; \mathbf{a}& =\mathbf{A} \mathbf{f} \\
\text{Synthesis:}\; \widehat{\mathbf{f}}& = \mathbf{S} \mathbf{a},
\end{align*}
where $\mathbf{a}$ is the vector of spherical harmonic coefficients and $\mathbf{f}$ and $\widehat{\mathbf{f}}$ are the vectors of data values at the HEALPix points. Using this notation, the iterative procedure in the HEALPix software can be written as
\begin{gather}
\begin{aligned}
    \mathbf{r}^{(k+1)} &= \mathbf{f} - \bS \mathbf{a}^{(k)}, \\
    \mathbf{a}^{(k+1)} &= \mathbf{a}^{(k)} + \bA \mathbf{r}^{(k+1)},
    \label{uzawa}
\end{aligned}
\end{gather}
where $\mathbf{r}$ is the residual vector and $\mathbf{a}^{(0)} = \bA \mathbf{f}$.  Substituting the first equation of \eqref{uzawa} into the last and using the fact that the analysis matrix satisfies $\bA = \frac{4\pi}{N}\bS^{*}$, gives the iteration
\begin{align}
    \mathbf{a}^{(k+1)} = \mathbf{a}^{(k)} + \frac{4\pi}{N}\bS^{*}(\mathbf{f} - \bS\mathbf{a}^{(k)}) = \frac{4\pi}{N}\bS^{*}\mathbf{f} + \left(\mathbf{I} - \frac{4\pi}{N}\bS^{*}\bS\right)\mathbf{a}^{(k)}.
    \label{richardson}
\end{align}
This is just stationary Richardson iteration (or Gradient Decent) with relaxation parameter $\frac{4\pi}{N}$ applied to the normal equations $\bS^{*}\bS\mathbf{a}= \bS^{*}\mathbf{f}$~\cite[pp.\ 44--45]{benzi2005}. 
Thus, the iterative procedure converges to the least squares solution to \eqref{healsph}, provided the spectral radius of $\mathbf{I}-\frac{4\pi}{N}\bS^{*}\bS$ is less than one.  The spectral radius also determines the convergence rate.  Since the HEALPIx points are equidistributed, we know that \eqref{healsph} converges to the integral \eqref{sphip} as $N\rightarrow\infty$ (in exact arithmetic)~\cite{HMS16}. Thus, the spectral radius of $\mathbf{I}-\frac{4\pi}{N}\bS^{*}\bS$ converges to $0$ as $N\rightarrow\infty$ and we expect the iteration \eqref{richardson} to converge more rapidly as $N$ increases. Table \ref{tab:table1} gives evidence of this result by displaying the spectral radius of $\mathbf{I}-\frac{4\pi}{N}\bS^{*}\bS$ for increasing values of $N$.

\begin{table}[h!]
  \begin{center}
    \begin{tabular}{|c|c|c|}
    \hline
     $N_{side}$ & \textbf{$N$} & \textbf{$\rho\left(\bI-\frac{4\pi}{N}\bS^{*} \bS\right)$} \\
      \hline
      \hline
     2 & 48 & 0.1986\\
    4 & 192 &0.0932\\
    8 & 768 & 0.0600\\
    16 & 3072 & 0.0475\\
    32 & 12288 &0.0421\\
    \hline
    \end{tabular}
     \caption{Spectral radius of the Richardson iteration matrix from \eqref{richardson} for different values of $N$.}
     \label{tab:table1}
  \end{center}
\end{table}

The default option in the HEALPix software sets the number of iterations of \eqref{richardson} to 3.  While this does improve the accuracy of computing the spherical harmonic coefficients, it adds to the cost, as each iteration requires doing an analysis and synthesis (\eqref{healsph} and \eqref{synthesis}) at a cost of $\mathcal{O}(\ell_{max}^3)$ operations each.  Since the solution converges to the least squares solution, one could improve the convergence of the Richardson iteration method by using algorithms like LSQR or conjugate gradient on the normal equations~\cite{PaigeSaunders}. 

\subsubsection{Pixel Weights and Ring Weights}
As an alternative to the iterative scheme, the HEALPix software also has the option of using  quadrature weights to improve the accuracy of the computation of the spherical harmonic coefficients. In this case, the equal weight quadrature approximation~\eqref{healsph} is generalized to 
\begin{align}
a_\ell^m=\sum_{i=1}^N w_i\overbar{Y}_\ell^m(\lambda_i,\theta_i)f(\lambda_i,\theta_i), \; 0 \leq \ell \leq \ell_{max}, -\ell \leq m \leq \ell,
\label{healsph_w}
\end{align}
where $w_i$ are the quadrature weights. There are two options for using quadrature weights.  The first is ``pixel weights'', which uses different weights for \textit{each HEALPix point}.  These weights are computed using the positive quadrature weight algorithm from~\cite{Potts}, which consists of solving a system involving a Gram matrix containing the spherical harmonics whose size is proportional to $N$~\cite{MartinRings}.  For large $N$, the weights are computed once and stored.  The second option is to use ``ring weights'', which use different weights for \textit{each ring} of the HEALPix point sets. The computation of the ring weights is done using similar ideas to the pixel weights, but a much smaller system has to be solved~\cite{MartinRings}. The new method introduced in this paper does not use quadrature weights directly, but instead computes the bivariate Fourier coefficients of the HEALPix data and then uses these to obtain the spherical harmonic coefficients. 

\section{HP2SPH}
\label{HP2SPH}
The algorithm presented here, named HP2SPH, introduces a new way to calculate the spherical harmonic coefficients of data sampled at the HEALPix points \eqref{healpts}. The outline for the algorithm is given in Algorithm \ref{alg:HP2SPH}, and each of the pieces are described below.
\begin{algorithm}
\caption{HP2SPH}
\label{alg:HP2SPH}
\begin{algorithmic}
\vspace{.2cm}
\STATE \textbf{Input:} Data sampled at the HEALPix point set of size $N$, $\{f_j\}$, $j=1,\ldots,N$.
\STATE \textbf{Output:} Approximate spherical harmonic coefficients, $\{a_{\ell}^{m}\}$, $0\leq \ell \leq \ell_{max}$, $-\ell \leq m \leq \ell$

\vspace{.3cm}
\STATE 1. Transform the data to a tensor product  latitude-longitude grid:
\STATE \hspace{2ex} (i) \hspace{0.7ex} Upsample the data in longitude from the northern ($\mathcal{N}$) and southern ($\mathcal{S}$) point sets using FFT
\STATE \hspace{2ex} (ii) \hspace{0.05ex} Shift (interpolate) the data from the equatorial ($\mathcal{E}$) point set so it is longitudinally aligned
\STATE 2. Compute the bivariate Fourier coefficients:
\STATE \hspace{2ex} (i) \hspace{0.7ex} Apply the DFS method
\STATE \hspace{2ex} (ii) \hspace{0.05ex} Apply the inverse NUFFT-II in latitude
\STATE \hspace{2ex} (iii) Apply the inverse FFT in longitude
\STATE 3. Obtain the spherical harmonic coefficients via the FSHT
\end{algorithmic}
\end{algorithm}

\subsection{Step 1: Transform the data to a tensor product  latitude-longitude grid}
\begin{figure}[ht]
	\begin{minipage}{0.5\textwidth}
		\centering
		\includegraphics[scale=.48]{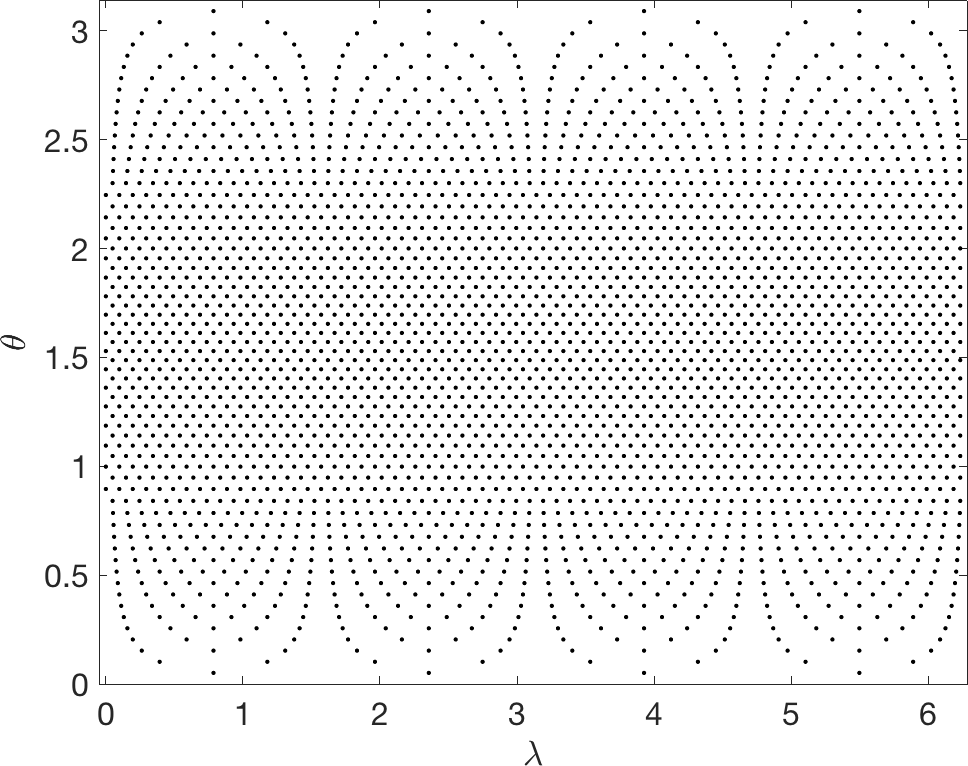}
		
		 \quad \;(a)
	\end{minipage}
	\begin{minipage}{0.5\textwidth}
		\centering
		\includegraphics[scale=.48]{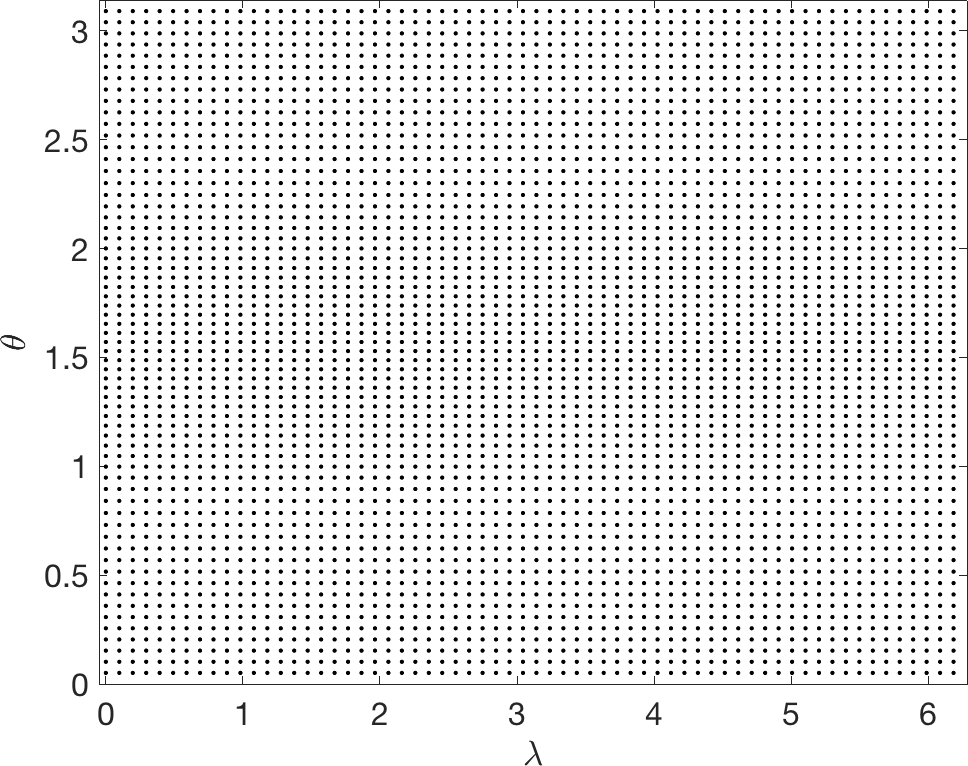}
		
		\quad \,(b)
	\end{minipage}
	\caption{(a) HEALPix points with $N_{side} = 16$ displayed in latitude and longitude and (b) the corresponding upsampled points.}
	\label{equiheal}
\end{figure}

As described in Section~\ref{HEALP}, the HEALPix grid has an unequal number of points on the rings in the northern ($\mathcal{N}$) and southern ($\mathcal{S}$) sets \eqref{healpts}, and the points on the rings in the equatorial ($\mathcal{E}$) set are shifted on every other ring. This structure leads to the pixels being misaligned in latitude. By upsampling the data on the northern and southern points in longitude so that there are an equivalent samples of the data on each ring and shifting the data at equatorial points in longitude, we can use fast algorithms to obtain the bivariate Fourier coefficients of the data as discussed in the next section. On the two polar point sets, we upsample the data using the trigonometric interpolant of the data on each ring of these sets to the non-shifted equally spaced longitude points on the equatorial rings, i.e.,\
\begin{equation}
\lambda_k = \frac{k}{2N_{side}}\pi, \; k=0,\dots,4N_{side}-1.
\label{upsamp}
\end{equation}
We also interpolate the data on the rings in the equatorial point set with shifted longitude points, to these $\lambda$ values. Figure~\ref{equiheal}(b) illustrates the upsampling procedure leading to a tensor product latitude-longitude grid of data of size $(4N_{side}-1) \times 4N_{side}$.

We describe the interpolation procedure here for the data in the northern point set $\mathcal{N}$. Consider the latitude values for the northern rings, $\theta_j=\arccos \left(1-\frac{j^2}{3N_{side}^2}\right), \; j=1,\dots,N_{side}$. We approximate the data in each ring using a trigonometric expansion of the form
\begin{equation}
 \text{$f(\lambda,\theta_j)=:$}f_j(\lambda)=\sum_{n=-2j}^{2j}c_n^{(j)}e^{{\rm i}n\lambda},
\label{trig_expansion}
\end{equation}
The coefficients in the expansion are determined by enforcing interpolation of the given data values 
\begin{equation*}
f\left(\frac{k+\frac{1}{2}}{2j}\pi,\theta_j\right),\; k=0,\dots,4j-1.
\end{equation*}
With the minor algebraic manipulation of \eqref{trig_expansion},
\begin{align*}
f_j\left(\frac{k+\frac{1}{2}}{2j}\pi\right)=\sum_{n=-2j}^{2j-1}c_n^{(j)}e^{{\rm i}n\frac{k+\frac{1}{2}}{2j}\pi}
=\sum_{n=-2j}^{2j-1}c_n^{(j)}e^{{\rm i}n\frac{\pi}{4j}}e^{{\rm i}n\frac{k}{2j}\pi}
=\sum_{n=-2j}^{2j-1}\tilde{c}_n^{(j)}e^{{\rm i}n\frac{k}{2j}\pi}, \, k=0,\dots, 4j-1,
\end{align*}
we see the interpolation conditions yield the system
\begin{equation}
\sum_{n=-2j}^{2j-1}\tilde{c}_n^{(j)}e^{{\rm i}n\frac{k}{2j}\pi}=f\left(\frac{k+\frac{1}{2}}{2j}\pi,\theta_j\right),\; k=0,\dots,4j-1,
\label{interp_cond}
\end{equation}
which can be computed using the inverse FFT. We can then obtain the Fourier coefficients $c_n^{(j)}$ in \eqref{trig_expansion} for the data at the non-shifted values through simple multiplication\footnote{Horner's rule (and Estrin's scheme for higher accuracy at small angles~\cite{estrin_60}) could also be used to implement the shift, avoiding loss of accuracy due to evaluation of high frequency complex exponentials.}. 
Finally, we pad the vector containing the coefficients $c_n^{(j)}$ with the appropriate number of zeros to get to a total of $4 N_{side}$, so that the expansion in longitude in each ring has the same number of Fourier coefficients.  The values of the interpolant on each ring can then be obtained at the upsampled values~\eqref{upsamp} by applying the FFT on these padded vectors. A similar procedure can be applied to the data on the rings in the southern point set $\mathcal{S}$. 

On the rings in the equatorial set $\mathcal{E}$ where the longitude values are shifted by $\pi(k+\frac12)/(2N_{side})$, we compute the Fourier coefficients of the data using \eqref{interp_cond} with $j = N_{side}$.  We then obtain the coefficients in \eqref{trig_expansion} at the non-shifted points from which the values can be computed using the FFT. No padding or upsampling is needed in this case.

\subsection{Step 2: Compute Bivariate Fourier Coefficients}
Bivariate Fourier analysis for data on the sphere requires the application of the DFS method to obtain periodicity of the data in latitude and to retain spherical symmetry. When we apply this method to the upsampled HEALPix data, there is an issue that the points in latitude are not equally-spaced, making the standard FFT unsuitable.  To bypass this issue we use an NUFFT. Both the DFS technique and NUFFT method we use are discussed below for completeness.

\subsubsection{Double Fourier Sphere (DFS) Method}
A natural approach to representing a function on the sphere is to use a latitude-longitude coordinate transform, defined by
\begin{equation}
x(\lambda,\theta)=\cos\lambda\sin\theta, \quad y(\lambda,\theta)=\sin\lambda\sin\theta \quad z(\lambda,\theta)=\cos\theta, \qquad (\lambda,\theta)\in [0,2\pi]\times[0,\pi],
\label{latlon}
 \end{equation}
which maps the sphere to a rectangular domain. While this transformation allows for performing computations with the function  $f(\lambda,\theta)=f(x(\lambda,\theta),y(\lambda,\theta),z(\lambda,\theta))$, it also introduces artificial boundaries at the north and south poles. In addition, the change of variables does not maintain the symmetry of functions on the sphere. Specifically, the transform described in~\eqref{latlon} does not preserve the periodicity in the latitude direction, which is necessary for bivariate Fourier analysis to be applicable and for results to make physical sense. These problems are solved by the DFS method.
\begin{figure}[ht]
 \centering
 \begin{minipage}{.49\textwidth} 
 \centering
   \begin{overpic}[width=.5\textwidth]{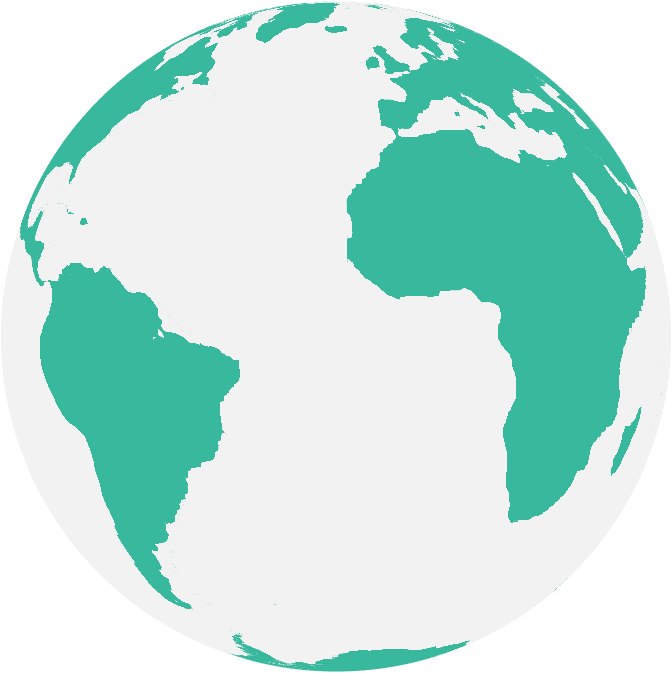}
   \end{overpic}
   
    \begin{overpic}[width=.95\textwidth,trim=30 50 30 60,clip]{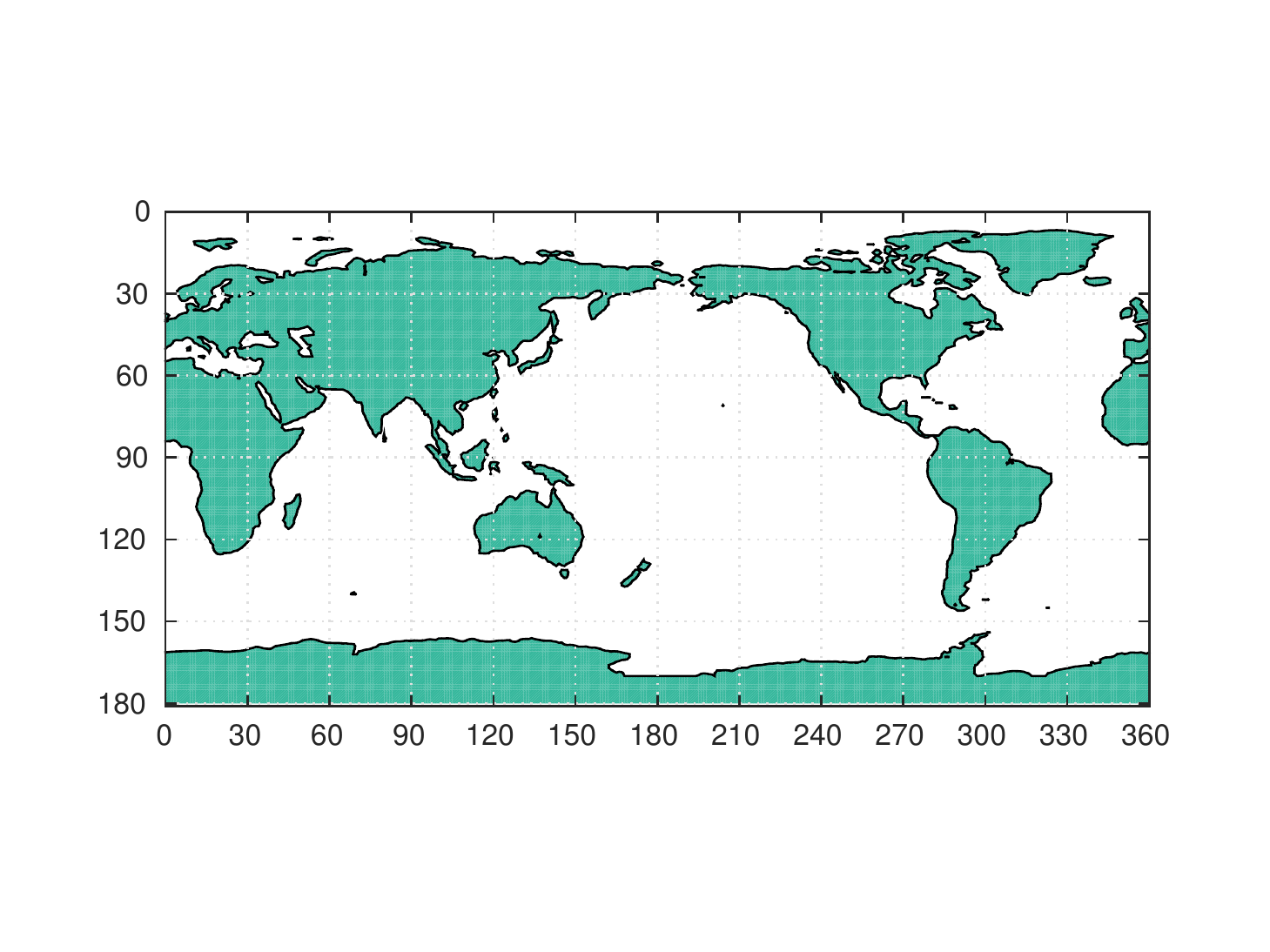}
    \put(0,105) {(a)}
    \put(0,60) {(b)}
    \put(50,0) {\footnotesize $180\lambda/\pi$}
    \put(-3,25) {\footnotesize \rotatebox{90}{$180\theta/\pi$}}
   \end{overpic}
 \end{minipage}
\begin{minipage}{.49\textwidth} 
   \begin{overpic}[width=\textwidth,trim=55 10 68 20,clip]{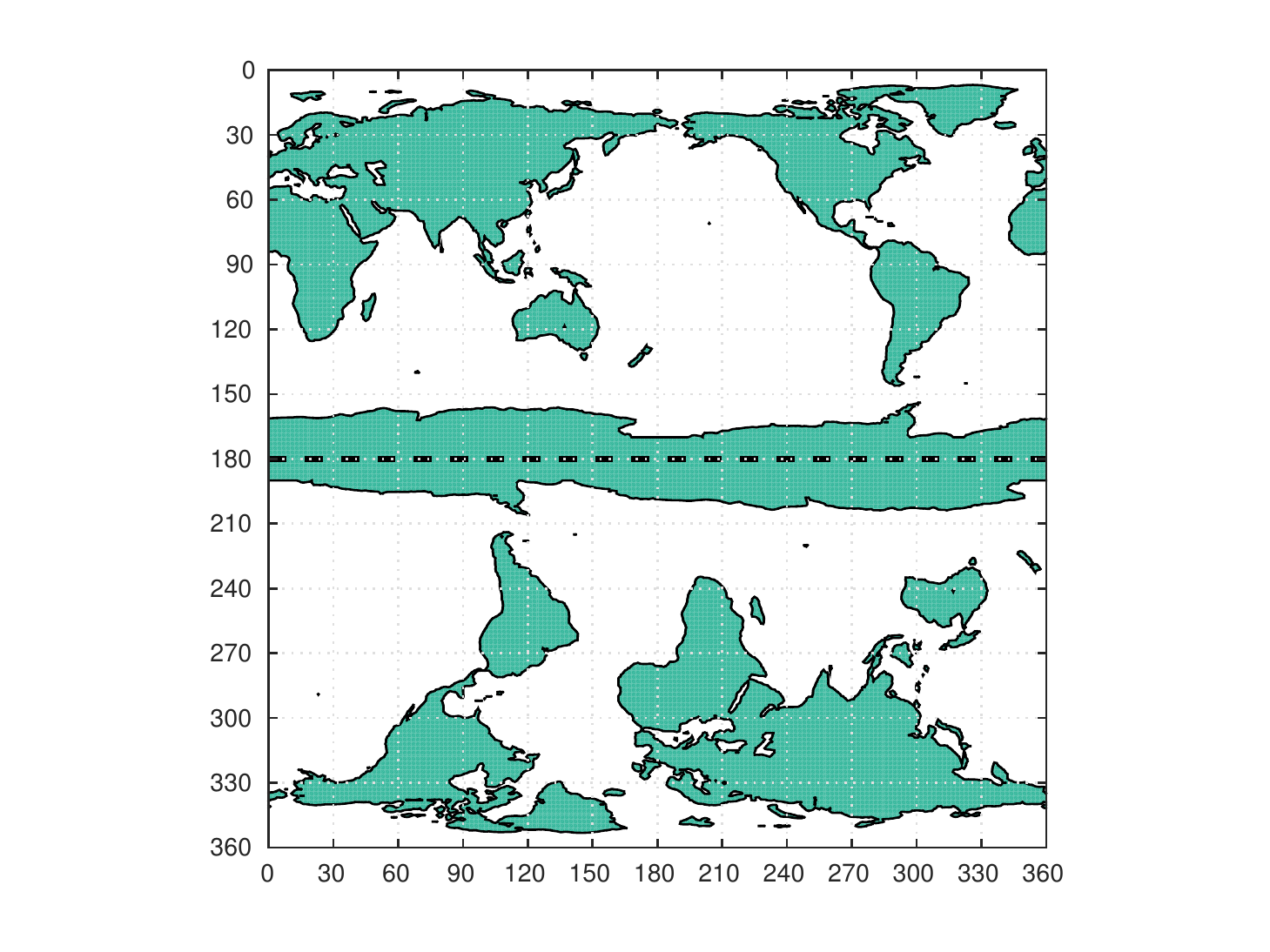}
  \put(0,95) {(c)}
  \put(55,0) {\footnotesize $180\lambda/\pi$}
  \put(0,45) {\footnotesize \rotatebox{90}{$180\theta/\pi$}}
 \end{overpic}
 \end{minipage}
 \caption{Illustration of the DFS method: (a) The surface of earth, (b) the surface mapped onto a latitude-longitude grid, and (c) the surface after applying the DFS method.~\cite{1sphere}}
 \label{atlas}
\end{figure}
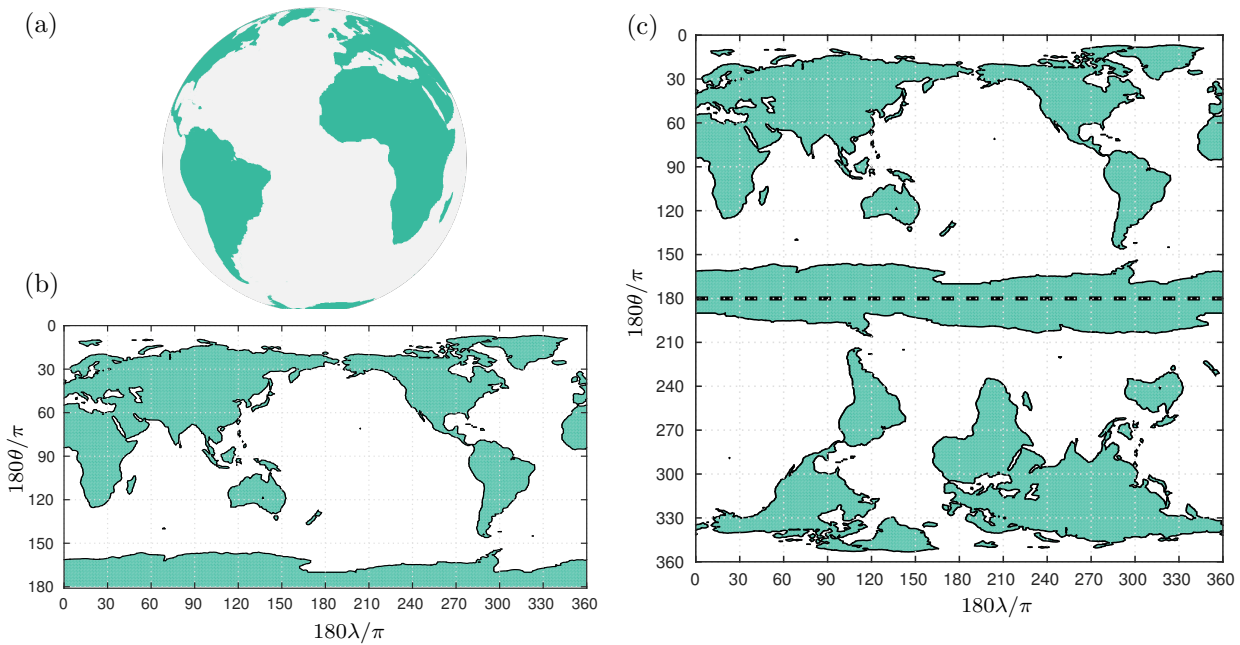

Originally introduced by Merilees in~\cite{DFS} (see also~\cite{1sphere}) the DFS method transforms a function on the sphere to a rectangular grid while maintaining bi-periodicity. This can be thought of as ``doubling up" the function $f(\lambda,\theta)$ to form a new function that preserves periodicity in both the latitude and longitude directions. Algebraically, this new function, $\widetilde{f}(\lambda,\theta)$, is defined on $[0,2\pi] \times [0,2\pi]$ as follows~\cite{1sphere}
\begin{equation}
\widetilde{f}(\lambda,\theta)=\begin{cases} 
      g(\lambda,\theta), & (\lambda,\theta)\in[0,\pi]\times [0,\pi],\\
      h(\lambda-\pi,\theta), & (\lambda,\theta)\in[\pi,2\pi]\times [0,\pi],\\
      h(\lambda,2\pi-\theta), & (\lambda,\theta)\in[0,\pi]\times [\pi,2\pi],\\
      g(\lambda-\pi,2\pi-\theta), & (\lambda,\theta)\in[\pi,2\pi]\times [\pi,2\pi],\\
   \end{cases}
   \label{DFS}
\end{equation}
where $g(\lambda,\theta)=f(\lambda,\theta)$ and $h(\lambda,\theta)=f(\lambda+\pi,\theta)$ for $(\lambda,\theta)\in[0,\pi]\times [0,\pi]$. Figure~\ref{atlas} illustrates the DFS method applied to the surface of the Earth, which shows the preservation of bi-periodicity in (c).  We note that the DFS method can also be easily applied to discrete data sampled at a tensor product latitude-longitude grid using \eqref{DFS}, which is what we do for the upsampled HEALPix data.  In this case, \eqref{DFS} corresponds to flipping and shifting the data matrix appropriately.

Once the DFS method is applied to a function on the sphere, it can be approximated using a 2D bivariate Fourier expansion:
\begin{equation}
\widetilde{f}(\lambda,\theta)\approx\sum_{j=-\lfloor \frac{m}{2} \rfloor}^{\lceil \frac{m}{2}\rceil-1}\,\sum_{k=-\lfloor \frac{n}{2} \rfloor}^{\lceil \frac{n}{2}\rceil-1} C_{jk} e^{{\rm i} j \theta}  e^{{\rm i} k \lambda},
\label{FFT2}
\end{equation}
where $m$ and $n$ represent the number of frequencies in (doubled-up) latitude and longitude, respectively.

Note that the HEALPix grid does not include points at the north and south poles.  When applying the DFS to the upsampled data from Step 1, this leads to a relatively large gap in the points in latitude over the poles compared to the other points, which can lead to issues with the inverse NUFFT (see below). To bypass this issue, we construct values at the two poles by using a weighted quadratic least squares fit~\cite{WLS} that combines the data from the three rings closest to each pole. 

\begin{remark}
The Fourier coefficients of the upsampled data in longitude are computed in Step 1.  These can be used directly in the DFS procedure by applying \eqref{DFS} in Fourier space in the $\lambda$ variable, which amounts to appending the (padded) coefficient matrix from Step 1 with a flipped version of itself with all odd wave numbers multiplied by $-1$. It then only remains to compute the Fourier coefficients in latitude $\theta$ to obtain the full bivariate Fourier coefficients.  This is the focus of the next step.
\end{remark}

\subsubsection{Nonuniform Fast Fourier Transform (NUFFT)}
The use of the nonuniform discrete Fourier transform (NUDFT) in many domain sciences has led to the development of algorithms for computing it efficiently. If these algorithms are quasi-optimal requiring $\mathcal{O}(n\log n)$, then they are referred to as a nonuniform fast Fourier transform (NUFFT). Given a vector $\mathbf{c}\in \mathbb{C}^{n}$, the one-dimensional NUDFT computes the vector $\mathbf{f}\in \mathbb{C}^{n}$ defined by
\begin{equation}
f_j=\sum_{k=0}^{n-1} c_k e^{-2\pi i x_j\omega_k}, \quad 0\leq j\leq n-1,
\label{DFT}
\end{equation}
where $x_j\in [0,1]$ are the samples and $\omega_k\in [0,n]$ are the frequencies. If the samples are equispaced ($x_j=j/n$) and the frequencies are integer ($\omega_k = k$), then the the transform is a uniform DFT, which can be computed by the FFT in $\mathcal{O}(n \log n)$ operations~\cite{fft2}. When either the samples are nonequispaced or the frequencies are noninteger, the FFT does not directly apply without some careful manipulation~\cite{fft1}.

Ruiz and Townsend~\cite{NUFFT} contributed to the collection of NUFFT algorithms by utilizing low rank matrix approximations to relate the NUDFT to the uniform DFT. There are three types of NUDFTs and inverse NUDFTs that they account for in their algorithm: NUDFT-I, which has uniform samples but noninteger frequencies; NUDFT-II, which has nonuniform samples and integer frequencies; NUDFT-III, which has both nonuniform samples and nonuniform frequencies~\cite{fft3}. Since our HP2SPH method only uses the one-dimensional inverse NUFFT of the second type, we discuss the NUFFT-II algorithm~\cite{NUFFT}.

Given Fourier coefficients, $\mathbf{c}\in\mathbb{C}^{n\times 1}$, the NUFFT-II attempts to approximate the vector
\begin{equation}
\mathbf{f}=\widehat{F}_2\mathbf{c},
\label{nufft2}
\end{equation}
to machine precision in quasi-optimal complexity.  Here $(\widehat{F}_2)_{jk}=e^{-2\pi i x_j k}$, $x_j$ are nonuniform samples (restricted to be in $[0,1]$), and $k$ are integer frequencies for $0\leq j,k \leq n-1$. Notice that the DFT matrix for uniform samples and integer frequencies is similarly $F_{jk}=e^{-2\pi ijk/n}$. The key to the NUFFT-II algorithm described in~\cite{NUFFT} is that if the samples are nearly equispaced, then $\widehat{F}_2$ can be related to the Hadamard product of $F$ and a low rank matrix. This means that given a rank $K$ approximation which relates $\widehat{F}_2$ and $F$, the NUFFT-II can then be computed using $K$ FFTs with $\mathcal{O}(Kn\log n)$ cost. In practice, machine (double) precision can be achieved with $K=14$~\cite{NUFFT}.

In the case of the inverse NUFFT-II, Ruiz and Townsend advocate for the use of the conjugate gradient (CG) method in order to solve the linear system $\widehat{F}_2\mathbf{c}=\mathbf{f}$ for $\mathbf{c}$. Since $\widehat{F}_2$ is not hermitian, the CG method is applied to the normal equations:
\begin{equation}
\widehat{F}_2^*\widehat{F}_2\mathbf{c}=\widehat{F}_2^*\mathbf{f},
\label{normal_eq}
\end{equation}
where $\widehat{F}_2^*\widehat{F}_2$ is simply a Toeplitz matrix. Therefore, the inverse NUFFT-II can be calculated using the CG method and a fast Toeplitz multiplication~\cite{toe} in $\mathcal{O}(R_{\scriptscriptstyle CG}n\log n)$ operations, where $R_{\scriptscriptstyle CG}$ is the number of CG iterations. The following suggestion is placed on the nonuniform function samples to avoid ill-conditioning in the system \eqref{normal_eq}~\cite{NUFFT}:
\begin{equation}
\left|x_j-\frac{j}{n}\right|\leq \frac{\gamma}{n}, \quad 0\leq j\leq n-1,
\label{equigrid}
\end{equation}
where $0\leq \gamma < 1/4$. When this condition is satisfied, it has been experimentally observed that $R_{\scriptscriptstyle CG}\leq 10$ for a large range of $n$.

For the method proposed in this paper, we apply the inverse NUFFT-II in latitude to the DFS upsampled HEALPix data from Step 2. Unfortunately, the HEALPix points in latitude direction do not meet the condition \eqref{equigrid}. To bypass this issue, we instead use a least squares solution to compute fewer coefficients at the higher wave numbers than what the given data may support. We describe this procedure below since it not discussed in~\cite{NUFFT}.

The inverse NUFFT-II method computes first column of the symmetric Toeplitz matrix $\widehat{F}_2^*\widehat{F}_2$ in \eqref{normal_eq} in the following manner:
 \begin{align*}
 \widehat{F}_2^*\widehat{F}_2 \mathbf{e}_1=\widehat{F}_2^*\boldsymbol{1} =(\boldsymbol{1}^{T}\widehat{F}_2)^*=(\widehat{F}_2^T\boldsymbol{1})^*.
\end{align*}
The last expression above can be computed efficiently by the NUFFT-I algorithm, since the NUDFT-I matrix is simply the transpose of the NUDFT-II matrix~\cite{NUFFT}.  To compute a least squares solution to~\eqref{nufft2} with fewer coefficients, we simply truncate the first column obtained from the above method to $m < n$ terms and form the resulting $m \times m$ Toeplitz matrix $\widehat{F}_2^*\widehat{F}_2$.  The right hand side for the least squares solution is obtained by similarly computing $\widehat{F}_2^*\mathbf{f}$ and truncating this to $m$ terms.

For the DFS upsampled HEALPix data from Step 2, there are $8N_{side}$ coefficients in latitude, but only $4N_{side}$ coefficients in longitude.  To keep the number of Fourier modes in both directions (nearly) the same, we choose $m=4N_{side}+1$ as the truncation parameter for the least squares solution for computing the Fourier coefficients in latitude.  This is also a convenient choice since the method in step three for converting bivariate Fourier coefficients of data on the sphere to spherical harmonic coefficients requires the number of coefficients in each direction is the same and an odd number (we explain how to convert the coefficients in longitude to $m=4N_{side}+1$ in the next section).

\begin{remark}
For problems where the data may contain noise (e.g., for the CMB application), there could be an issue with this noise being amplified in steps 1 and 2.  For step 1, we should not expect any additional noise to be introduced, since we are simply computing the Fourier coefficients on each ring using the original data and then shifting the coefficients and padding them with zeros. Step 2 has two areas where there could be an issue with noisy data.  The first is in constructing values at the poles and the second is in the application of the NUFFT in latitude. However, both of these steps apply a least squares procedure, which provides some smoothing.  In our tests on CMB data, we did not observe any amplification of noise that was present in the data.  The HP2SPH method has a further benefit of using a backward stable algorithm for computing the spherical harmonic coefficients (as discussed next), which ensures that the resulting uncertainty in the spherical harmonic coefficients has only a low algebraic growth with respect to degree and is always proportional to the norm of the noise in the data.
\end{remark}

\subsection{Step 3: Obtain the spherical harmonic coefficients via the fast spherical harmonic transform (FSHT)}
In~\cite{FastSPH}, Slevinsky derives a fast, backward stable method for the transformation between spherical harmonic expansions and their bivariate Fourier series (given in~\eqref{FFT2}) by viewing it as a change of basis. This relation is defined as a connection problem, and the matrices that arise in the present case are well-conditioned, making them ideal for fast computations. Slevinsky describes the change of basis in two steps: converting expansions in normalized associated Legendre functions to those of only order zero and one, and then re-expressing these in trigonometric form. In other words, it uses spherical harmonic expansions of order zero and one as intermediate expressions between higher-order spherical harmonics and their corresponding bivariate Fourier coefficients. 

The first step of the algorithm takes advantage of the fact that the matrix of connection coefficients between the associated Legendre functions of all orders and those of order zero and one can be represented by a product of Givens rotation matrices. This enables the use of the butterfly algorithm, which can be thought of as an abstraction of the algebraic properties of fast Fourier transforms. The term butterfly was introduced in~\cite{butterfly1}, where it was used for analyzing scattering from electrically large surfaces, and then further developed in~\cite{butterfly2} for use in special function transforms. Slevinsky uses the butterfly algorithm to perform a factorization of the well-conditioned matrices of connection coefficients.

The second step of this method exploits the hierarchical decompositions of the connection coefficient matrices between the associated Legendre functions of order zero and one to the Chebyshev polynomials of the first and second kind, respectively. This step quickly computes the fast orthogonal polynomial transforms using an adaptation of the Fast Multipole Method~\cite{FMM} and low-rank matrix approximations. The total pre-computation time for both steps is $\mathcal{O}(\ell_{max}^3\log \ell_{max})$, and execution time is asymptotically optimal with $\mathcal{O}(\ell_{max}^2\log^2 \ell_{max})$ operations. This FSHT is implemented in Julia with the package FastTransforms~\cite{fasttrans} (as are the NUFFT methods from~\cite{NUFFT} used in Step 2). 

The FSHT in FastTransforms assumes the input function has a bivariate Fourier expansion of the form
 \begin{equation}
  \widetilde{f}(\lambda,\theta)=\sum_{j=0}^p\sum_{k=-p}^{p} g_{j}^k \frac{e^{{\rm i} k \lambda}}{\sqrt{2\pi}} \left\{\begin{array}{c}\cos j\theta, \quad \quad \; k\, \text{even}\\ \sin(j+1)\theta, \, k\, \text{odd}\;\end{array}\right\}.
  \label{bivariateFourierSph}
  \end{equation}
Any function on the sphere is required to have these even/odd conditions on its bivariate Fourier coefficients~\cite{DFS}.  At the end of step 2 we have obtained the bivariate Fourier expansion of the data of the form
\begin{align}
\widetilde{f}(\lambda,\theta)=\sum_{j=-p}^{p}\,\sum_{k=-p}^{p-1} C_{jk} e^{{\rm i} j \theta}  e^{{\rm i} k \lambda},
\label{bivariateFourier}
\end{align}
where $p = N_{side}/2$. 
Since we are dealing with real-valued data, we can expand Fourier coefficients array in $\lambda$ to an odd number of points.  The expanded array is defined by
\begin{align*}
X_{j,k} = 
\begin{cases}
    C_{j,k} & \text{if $-p+1 \leq k \leq p-1$} \\
    \frac{1}{2} C_{j,-p} & \text{if $k=\pm p$} \\
\end{cases},\;
-p\leq j,k \leq p.
\end{align*}
Using the array $X$, we can write \eqref{bivariateFourier} as
\begin{align*}
 \widetilde{f}(\lambda,\theta)&=\sum_{k=-p}^{p} e^{{\rm i} k \lambda} \sum_{j=0}^{p}((X_{jk}+\overbar{X}_{-j k})\cos(j\theta))+\frac{1}{{\rm i}}(X_{jk}-\overbar{X}_{-jk})\sin(j\theta))\\
&= \sum_{j=0}^p\sum_{k=-p}^{p} e^{{\rm i} k \lambda}\left\{\begin{array}{c}((X_{jk}+\overbar{X}_{-jk})\cos(j\theta),\; \text{$k$ even}\\ ((X_{jk}-\overbar{X}_{-jk})\sin(j\theta),\; \text{$k$ odd}\end{array}\right\},
\end{align*}
from which we can obtain the coefficients $g_j^k$ in \eqref{bivariateFourierSph}. 

The FSHT software takes bivariate Fourier coefficients $g_j^k$ as input in an array organized as follows:
\begin{equation*}
\begin{pmatrix}
g_0^0 & g_0^{-1} & g_0^1 & \cdots & g_0^{-p} & g_0^{p}\\
g_1^0 & g_1^{-1} & g_1^1 & \cdots & g_1^{-p} & g_1^{p}\\
\vdots & \vdots & \vdots & \ddots & \vdots & \vdots\\
g_{p-1}^0 & g_{p-1}^{-1} & g_{p-1}^1& \cdots & g_{p-1}^{-p} & g_{p-1}^{p}\\
g_p^0 & 0 & 0 & \cdots & g_p^{-p} & g_p^{p}\\
\end{pmatrix}.
\end{equation*}
The output of the software is the approximate spherical harmonic coefficients of the data arranged in an array of the form
\begin{equation*}
\begin{pmatrix}
a_0^0 & a_1^{-1} & a_1^1 & a_2^{-2} & a_2^2 & \cdots & a_p^{-p} & a_p^p\\
a_1^0 & a_2^{-1} & a_2^1 & a_3^{-2} & a_3^2 & \cdots & 0 & 0\\
\vdots & \vdots & \vdots &  \vdots &  \vdots & \ddots & \vdots & \vdots\\
a_{p-2}^0 & a_{p-1}^{-1} & a_{p-1}^1 & a_p^{-2} & a_p^2 &  & \vdots & \vdots\\
a_{p-1}^0 & a_p^{-1} & a_p^1 & 0 & 0 & \cdots & 0 & 0\\
a_p^0 & 0 & 0 & 0 & 0 & \cdots & 0 & 0\\
\end{pmatrix}.
\end{equation*}
The angular power spectrum \eqref{APS} can then be computed from this array.

\section{Numerical Results}
\label{NumRes}
In this section we present a few numerical tests comparing the spherical harmonics and angular power spectrum \eqref{APS} computed by our new method HP2SPH to the values computed by the HEALPix software employing the iterative scheme~\eqref{richardson}, pixel weights, and ring weights~\eqref{healsph_w}.  The first test compares the rate at which the two methods converge to the spherical harmonic coefficients by applying them to deterministic (i.e.\ non-noisy) functions sampled at the HEALPix points with known coefficients. The second test compares the accuracy of the methods using deterministic functions, which have a known power spectrum. In the third test, we compare the methods after calculating the angular power spectrum for a real CMB data map, which contains noise.
\subsection{Convergence of Spherical Harmonic Coefficients}
\begin{figure}[ht]
\centering
\begin{tabular}{cc}
		\includegraphics[width=0.47\textwidth]{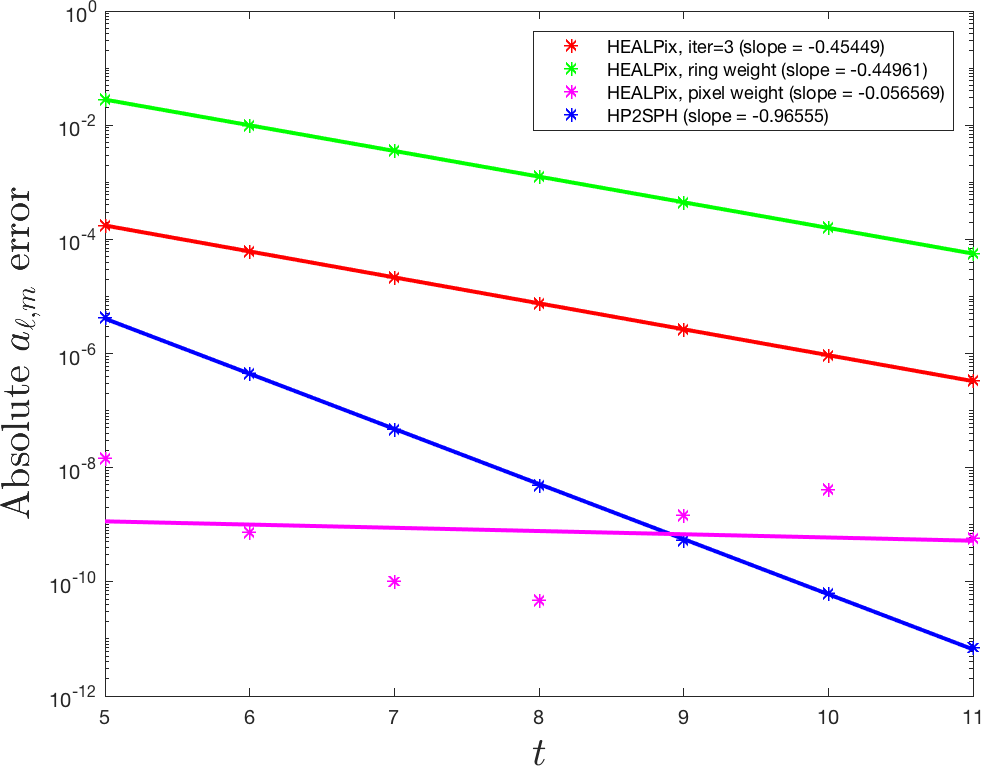} & 
		\includegraphics[width=0.47\textwidth]{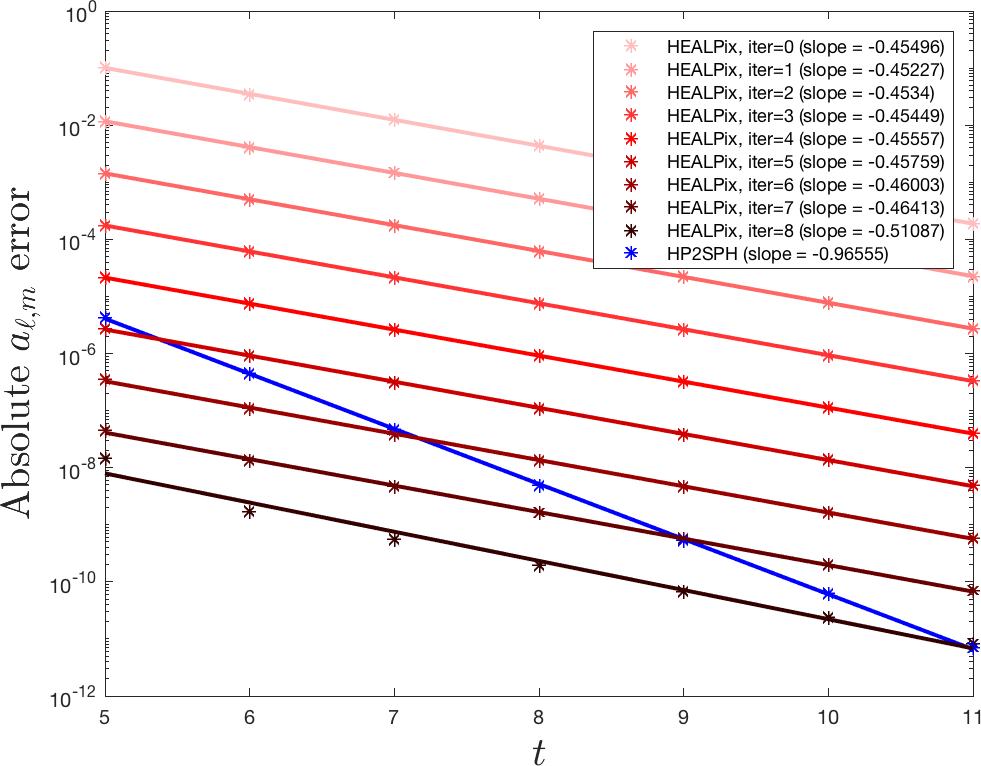} \\
		(a) & (b) \\
\end{tabular}		
	\caption{Maximum absolute errors as a function of $t$ for the computed spherical harmonic coefficients of~\eqref{f_plot} using HP2SPH and (a) HEALPix (3 iterative refinement steps), pixel weights, ring weights and (b) HEALPix with increasing iterative steps. The lines in the figure are the lines of best fit to the data and the convergence rates are determined from the slope of this line (as displayed in the plot legends).}
	\label{converge_plots}
\end{figure}

We choose the test function
\begin{equation}
f(\lambda,\theta)=\sum_{j=1}^3 c_j (2 - 2\mathbf{x}(\lambda,\theta) \cdot \mathbf{x}(\lambda_j,\theta_j))^{3/2},
\label{f_plot}
\end{equation}
where $\mathbf{x}(\lambda,\theta)=[x(\lambda,\theta)\; y(\lambda,\theta)\;z(\lambda,\theta)]$ from \eqref{latlon} and the parameters, which we picked randomly, are given by
\begin{gather*}
\{c_1,c_2,c_3\} = \{5,-3,8\}, \\
\{\lambda_1,\lambda_2,\lambda_3\} =  \{0.891498158152027,2.650004294134628,5.753735997130328\}, \\
\{\theta_1,\theta_2,\theta_3\} = \{1.232217523107963,2.059244524372349,0.537798840821172\}.
\end{gather*}
The function $(2 - 2\mathbf{x}(\lambda,\theta) \cdot \mathbf{x}(\lambda_{\rm c},\theta_{\rm c}))^{3/2}$ is a called a \emph{potential spline} of order $3/2$ centered at $\mathbf{x}(\lambda_{\rm c},\theta_{\rm c})$ and its exact spherical harmonic coefficients are given by~\cite{BaHu}
\begin{align}
    \widetilde{a}_\ell^m=\frac{18\pi}{(\ell+5/2)(\ell+3/2)(\ell+1/2)(\ell-1/2)(\ell-3/2)}Y_\ell^m(\lambda_{\rm c},\theta_{\rm c}).
    \label{alm}
\end{align}
These values are used to compare the convergence rates of the methods to the exact spherical harmonic coefficients of $f$. We do this by plotting in Figure \ref{converge_plots} the maximum absolute errors of the coefficients against the parameter $t$, which is used to determine the grid resolution parameter ($N_{side} = 2^t$).
Note that the spherical harmonic coefficients of the CMB decay at a rate between $\mathcal{O}(\ell^{-2})$ and $\mathcal{O}(\ell^{-3})$~\cite{low_ell}, which is slower than the decay rate of the coefficients of the test function~\eqref{f_plot} (since, for all $-\ell \leq m\leq \ell$, $|Y_\ell^m(\lambda_{\rm c},\theta_{\rm c})| \leq \sqrt{(2\ell + 1)/4\pi}$, and the remaining terms in \eqref{alm} decay at a rate of $\mathcal{O}(\ell^{-5})$).  This means that the test function has more smoothness than an actual CMB data set.

Figure~\ref{converge_plots}(a) compares the four methods and shows that the HP2SPH method converges to the spherical harmonic coefficients of~\eqref{f_plot} at a rate at least twice as fast as any of the HEALPix methods. Although the consecutive iterative refinement steps used in the HEALPix method produce progressively better errors, Figure~\ref{converge_plots}(b) illustrates that this does not improve the convergence rate (as discussed in Section~\ref{HEALsoft}). It is also important to note that after 8 iterative steps, there are no further improvements in the accuracy, indicating the algorithm has nearly converged to the least squares solution to \eqref{healsph}. The results the HEALPix method with pixel weights look pretty erratic with convergence achieved to 8-10 digits around $t=7$, but no further reductions. This could be because of potential errors in the computed quadrature weights.

Next, we test how the convergence rates of the methods are affected by high frequencies.  To do this we add 15 spherical harmonics to~\eqref{f_plot} with the following randomly chosen degrees and orders:
\begin{table}[h]
\centering
\begin{tabular}{|c|c|c|c|c|c|c|c|c|c|c|c|c|c|c|c|}
\hline
$\ell$ & 176 & 190 & 191 & 230 & 248 & 283 & 292 & 303 & 326 & 366 & 388 & 404 & 421 & 446 & 448   \\ \hline
$m$ & 56 & 81 & 124 & 40 & 155 & 274 & 27 & 145 & 55 & 343 & 200 & 78 & 420 & 284 & 234     \\ \hline
\end{tabular}
\end{table}
\begin{figure}[h!]
		\centering
		\includegraphics[scale=.5]{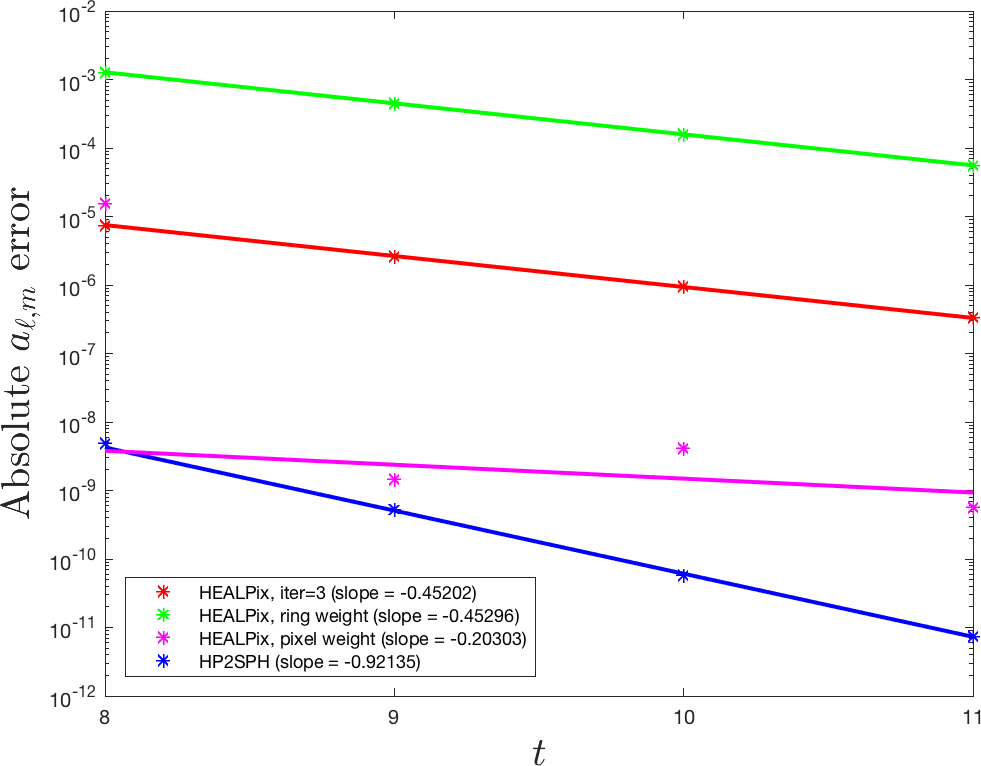}
		\caption{Maximum absolute errors as a function of $t$ for the computed spherical harmonic coefficients of \eqref{f_plot} augmented with spherical harmonics using HP2SPH, HEALPix with 3 iterative refinement steps, pixel weights, and ring weights.}
		\label{sph_conv}
\end{figure}
\noindent The results from this test are displayed in Figure~\ref{sph_conv}. We note that in order to ensure the calculation of asymptotic convergence rates for all methods, we excluded the errors for $t=8$ in the lines of best fit. This test shows that the convergence rates of the new method as well as the iterative HEALPix scheme and the HEALPix method with ring weights are not affected by the addition of high degree spherical harmonic terms to~\eqref{f_plot}. The HEALPix method with pixel weights shows a similar erratic behavior to Figure~\ref{converge_plots}(a), with no steady reductions in the errors after $t=9$.  The new HP2SPH method has the lowest errors of all the methods.

\subsection{Errors in the Angular Power Spectrum}
In this test, we first explore the accuracy of all the methods for computing the angular power spectrum of~\eqref{f_plot}. These results are compared to the exact spectrum, which is calculated using the exact spherical harmonic coefficients~\eqref{alm} in~\eqref{APS}. 
\begin{figure}[ht]
\centering
\begin{tabular}{cc}
		\includegraphics[width=0.44\textwidth]{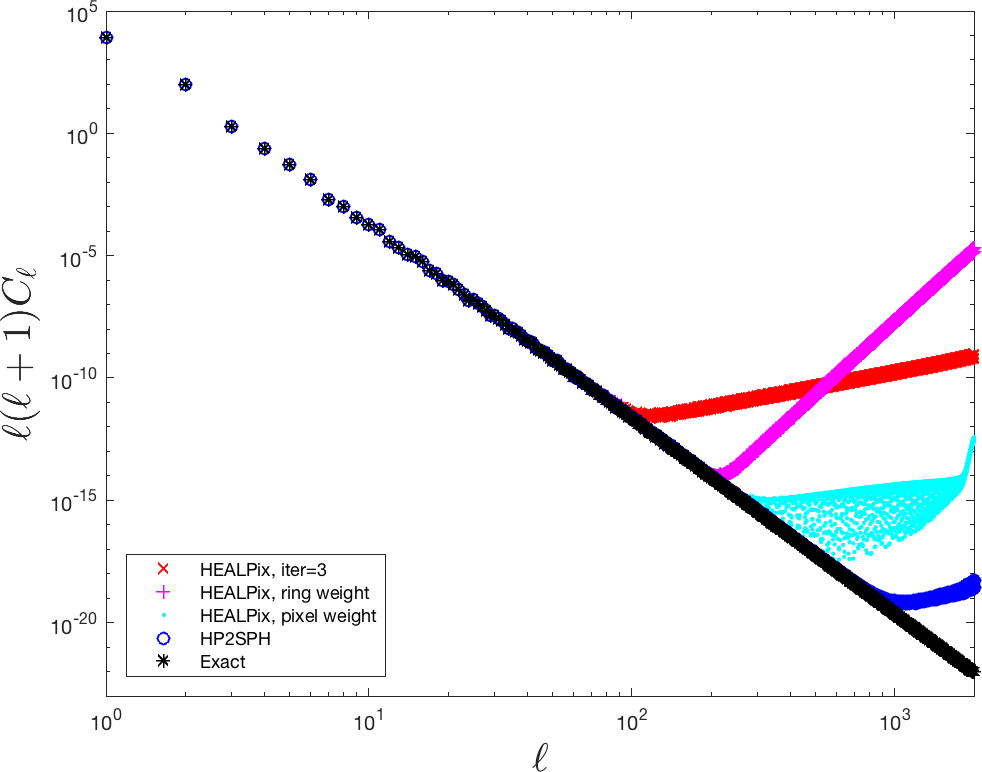} & 
		\includegraphics[width=0.44\textwidth]{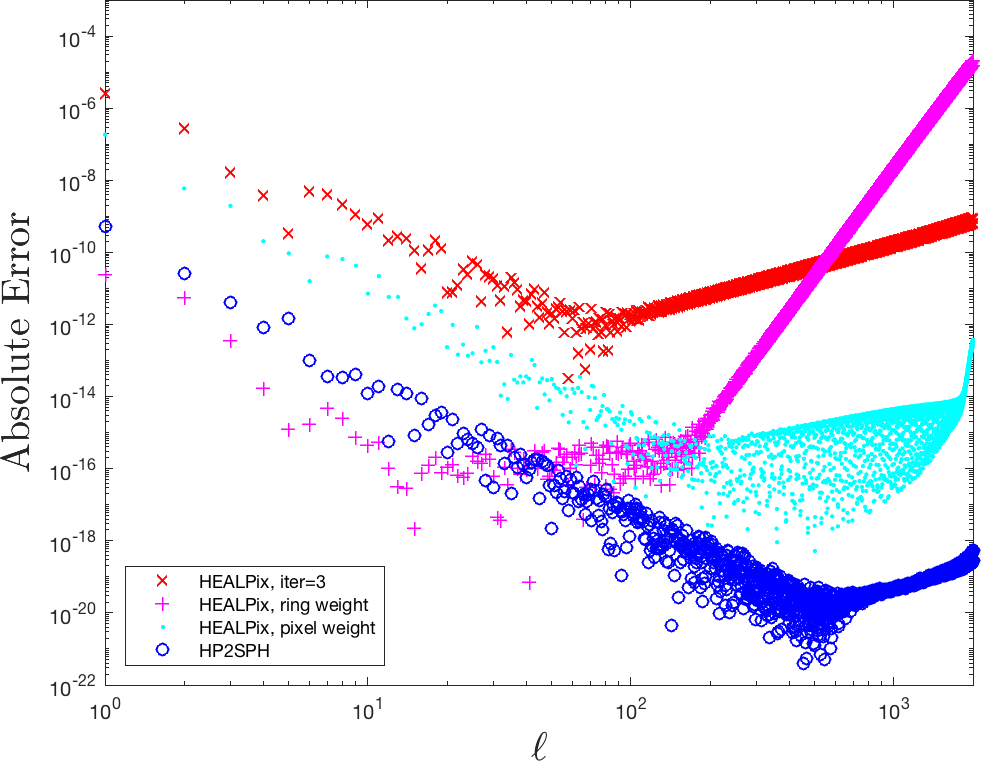} \\
		(a) & (b) \\
\end{tabular}		
	\caption{(a) Scaled angular power spectrum of \eqref{f_plot} as a function of degree $\ell$ computed by the HEALPix software with 3 iterative refinement steps, the HP2SPH method, the HEALPix method with ring weights, and the HEALPix method with pixel weights. The exact power spectrum is given by the black $\circ$'s. Here $N_{side}=2^{10}$, which is $N=12,\;582,\;912$ total points. (b) Absolute errors in the (scaled) angular power spectrum of the results from (a) as a function of degree $\ell$.}
	\label{f_aps}
\end{figure}
The angular power spectrum of the four methods are displayed in Figure~\ref{f_aps}(a).  We see that the algorithms produce similar results for lower degrees $\ell$, but the HEALPix method with the iterative refinement scheme~\eqref{richardson} diverges for degrees greater than approximately $\ell=100$, whereas the ring weight and pixel weight quadrature~\eqref{healsph_w} results diverge for degrees greater than $\ell=200$. To compare the methods more directly, Figure~\ref{f_aps}(b) plots the absolute errors in the angular power spectrum for each degree $\ell$. While the HEALPix method with ring weights performs the best out of all of the HEALPix methods for $\ell \leq 50$, the errors increase rapidly for higher $\ell$. Conversely, the HEALPix method with pixel weights does not perform as well for smaller $\ell$, yet it performs better than the other HEALPix methods for larger $\ell$. The HP2SPH method offers comparable results for low $\ell$ as the pixel weights method while still maintaining accuracy for high $\ell$.

Similar to Test 1, we test the accuracy of the methods when computing the power spectrum of data with high frequencies, as occur in real CMB data maps. As before, we add several spherical harmonic functions of high degree to the function~\eqref{f_plot}. The new test function appends the following degrees and orders:

\begin{table}[h]
\begin{tabular}{|c|c|c|c|c|c|c|c|c|c|c|c|c|c|c|c|}
\hline
$\ell$ & 589 & 633 & 636 & 766 & 829 & 943 & 974 & 1009 & 1085 & 1219 & 1294 & 1346 & 1404 & 1485 & 1493   \\ \hline
$m$ & 188 & 269 & 414 & 134 & 517 & 912 & 93 & 483 & 183 & 1143 & 667 & 259 & 1400 & 946 & 779   \\ \hline
\end{tabular}
\end{table}

\noindent The power spectrum of this function is the same as \eqref{f_plot}, but with the value at each degree $\ell$ of appended spherical harmonics increased by $\frac{1}{2\ell +1}$.

\begin{figure}[h!]
\centering
\begin{tabular}{cc}
		\includegraphics[width=0.44\textwidth]{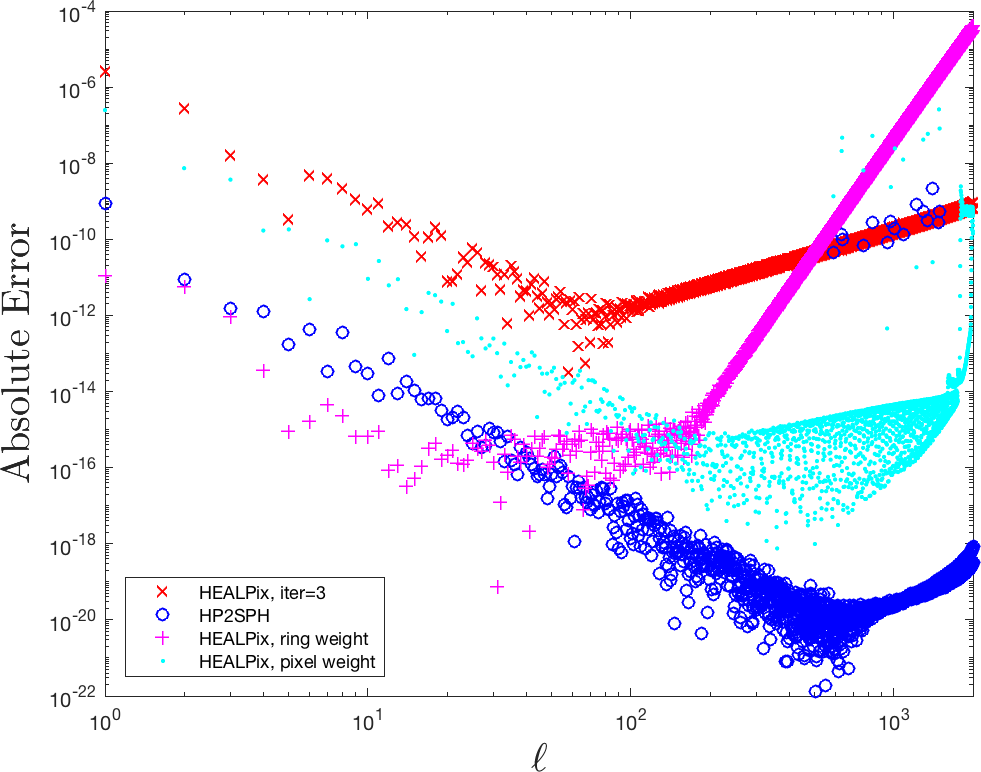} & 
		\includegraphics[width=0.44\textwidth]{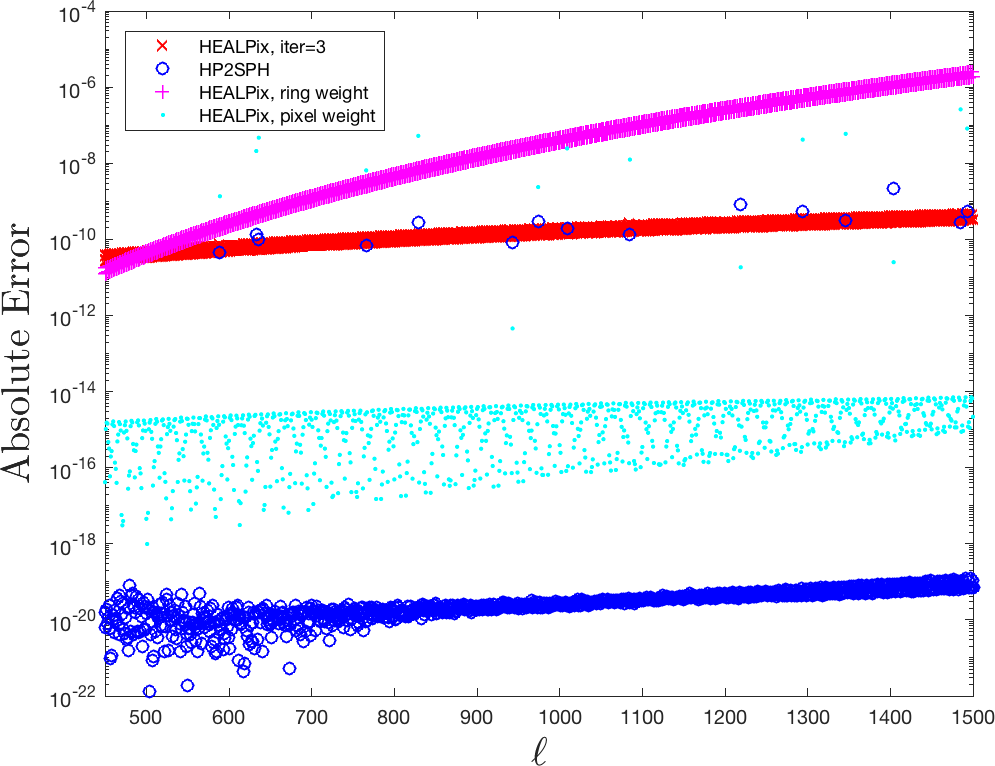} \\
		(a) & (b)
\end{tabular}		
	\caption{Absolute errors in the (scaled) angular power spectrum of \eqref{f_plot} augmented with high-degree spherical harmonics computed by the HEALPix software with 3 iterative refinement steps, the HP2SPH method, the HEALPix method with ring weights, and the HEALPix method with pixel weights as a function of degree $\ell$. (a) Displays the errors for degrees $\ell=1,\dots,2000$, while (b) displays the errors only for $\ell=450,\dots,1500$ to better show the how good the methods are at recovering the spectrum at the degrees of the appended spherical harmonics. Here $N_{side}=2^{10}$, which is $N=12,\;582,\;912$ total points.}
	\label{sph_test}
\end{figure}

Figure~\ref{sph_test}(a) displays the errors in the angular power spectrum of this new function for all the methods over all $\ell$, while Figure~\ref{sph_test}(b) displays the errors only over the range of $\ell$ that were appended to the base function.  We again see that the ring weights provide the highest accuracy for low $\ell$, but then the errors increase rapidly, while the HEALPix iterative method has the largest errors for low $\ell$, but some of the smallest errors for the $\ell$ corresponding to the appended spherical harmonics. In contrast, the new HP2SPH method provides small errors over the entire angular power spectrum, clearly showing its benefits over the HEALPix methods for determining the angular power spectrum of deterministic functions on the sphere.

\subsection{Application to Real CMB Map}
\begin{figure}[h!]
\centering
\begin{tabular}{cc}
		\includegraphics[width=0.45\textwidth]{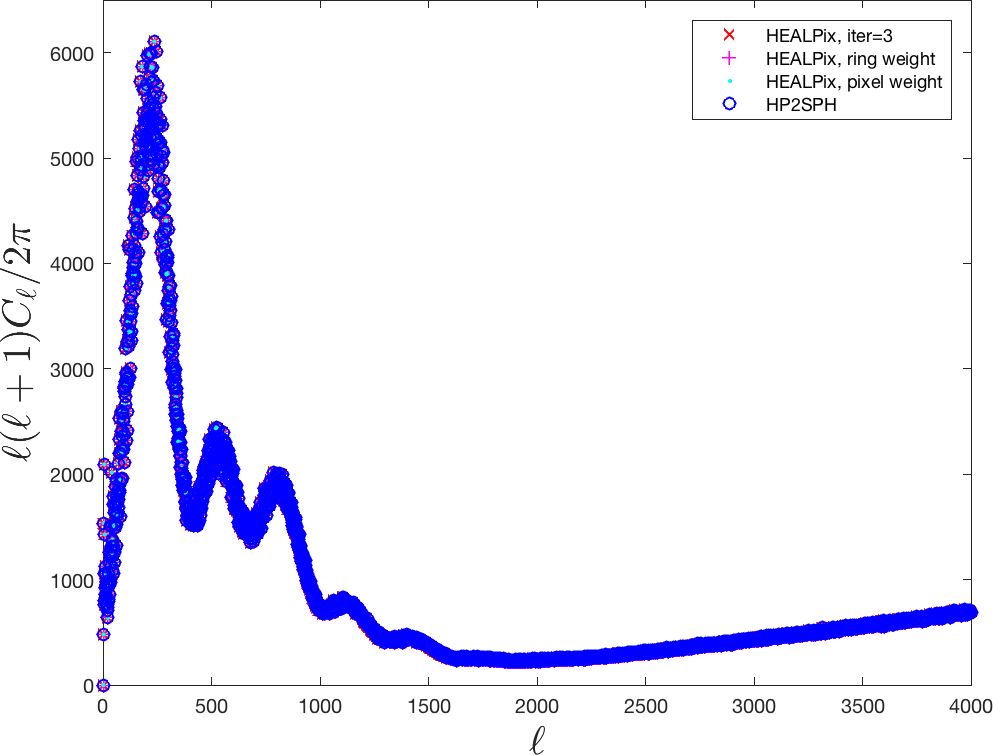} & 
		\includegraphics[width=0.45\textwidth]{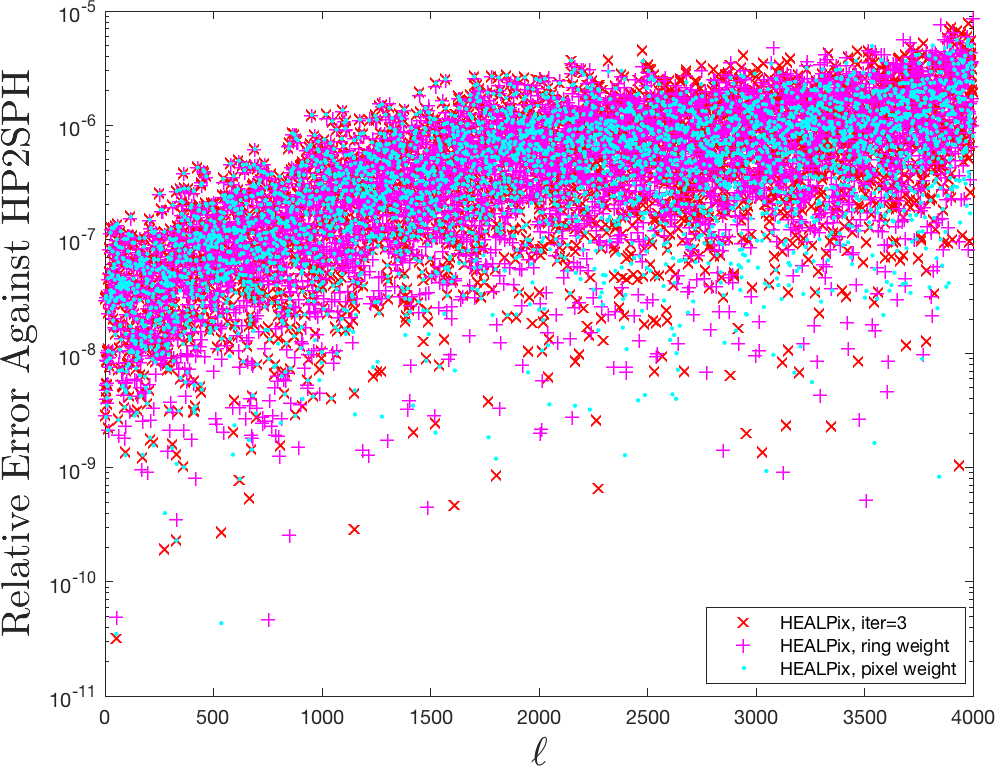} \\
		(a) & (b)
\end{tabular}		
	\caption{(Scaled) angular power spectrum of the CMB data map displayed in Figure \ref{CMB} (a) with $N_{side}=2^{11}$ for the four methods discussed in the paper (left), and the relative errors of the HEALPix software methods against the HP2SPH method (right).}
	\label{1024_map}
\end{figure}

Our final numerical test compares the methods on the real CMB map shown in Figure~\ref{CMB}. Figure \ref{1024_map} (a) shows the angular power spectrum for this map computed with the methods, while (b) shows the errors in the three HEALPix methods compared to the new HP2SPH method.  We see from the figure that the new method (in blue) produces visually the same results as the HEALPix methods (red, magenta, and cyan), indicating that the new method is not anymore susceptible to noise than the HEALPix methods.

\section{Conclusions and Remarks}
\label{Conc}
We have presented a new method, HP2SPH, for performing discrete spherical harmonic analysis on data collected using the HEALPix scheme. The method utilizes the FFT, NUFFT, and the FSHT to compute the spherical harmonic transform in near optimal computational complexity ($\mathcal{O}(N\log^2 N)$ complexity for $N$ total HEALPix points). Several numerical tests were presented to demonstrate the improved convergence and accuracy of the new method relative to the various HEALPix approaches for problems involving synthetic data with no noise, except that introduced by roundoff errors. In the case of a real CMB map with additional types of noise, the power spectra computed by the methods show good agreement. The new HP2SPH benefits further from the backward stability of the FSHT, which ensures that the resulting uncertainty in the spherical harmonic coefficients has only a low algebraic growth with respect to degree and is always proportional to the norm of the  noise in the  data.  We anticipate the new method will be applicable to the many other areas where the HEALPIx scheme is used and is naturally generalizable to other equal-area isolatitudinal sampling schemes for the sphere.

For our next steps, we will work to optimize the implementation of the method, which is currently in Julia, to improve its actual run-time.  This will include transcribing our code into a lower-level language like C; efforts in this direction are already underway for the FSHT~\cite{fasttransC}. In addition to this, we will include the ability to perform Fourier synthesis on a CMB map, i.e.\ given an angular power spectrum, we will return the corresponding CMB map values. For this purpose, our method has another advantage over HEALPix in that we will have the bivariate Fourier coefficients, which will simply make the synthesis an application of the FFT and NUFFT. Finally, we plan to add functionality for analyzing the polarization of CMB maps. 

\section*{Acknowledgements}
We are grateful to Ann Almgren for her comments and suggestions as well as Richard Mikael Slevinsky for assisting us with using the FastTransforms package. We also thank the entire HEALPix team, especially Kris G\'orski and Martin Reinecke, for their feedback and assistance on this work. Finally, we thank the anonymous reviewers whose comments helped improve the paper. The first author's work was supported by the SMART Scholarship, which is funded by The Under Secretary of Defense-Research and Engineering, National Defense Education Program / BA-1, Basic Research. The second author's work was supported by National Science Foundation grant CCF 1717556. 

\appendix

\section{Spherical Harmonic Conventions}\label{app:sph}
We denote a scalar spherical harmonic of degree $\ell\geq 0$ and order $-\ell\leq m\leq \ell$ as $Y_\ell^m(\lambda,\theta)$, where $\lambda$ is the azimuth angle and $\theta$ is the zenith angle. We define these functions as
\begin{equation}
	Y_\ell^m(\lambda,\theta) = \sqrt{\frac{2\ell+1}{4\pi}}\sqrt{\frac{(\ell-m)!}{(\ell+m)!}}P_\ell^m(\cos \theta)e^{i m \lambda},\quad m=0,1,\dots,\ell, \\ 
\label{ssph}
\end{equation}
where $Y_\ell^m=(-1)^mY_\ell^{|m|}$ for $m< 0$ and $P_\ell^m(\cos \theta)$ are the associate Legendre functions. As eigenfunctions of the Laplace-Beltrami operator, spherical harmonics are the natural basis for square integrable functions on the sphere~\cite{AtkinsonHan2012}. In other words, any $L^2$-integrable function $f$ on the sphere can be uniquely represented as
\begin{equation*}
f(\lambda,\theta) = \sum_{\ell=0}^\infty \sum_{m=-\ell}^\ell \widetilde{a}_\ell^mY_\ell^m (\lambda,\theta),
\end{equation*} 
where the spherical harmonic coefficients, $\widetilde{a}_\ell^m$, are found using the usual $L^2$-inner product for scalar functions on the sphere:
 \begin{equation}
 \widetilde{a}_\ell^m=\langle f,Y_\ell^m \rangle = \int_0^{2\pi}\int_0^\pi f(\lambda,\theta)\overbar{Y}_\ell^m(\lambda,\theta)\sin \theta d\theta d\lambda.
 \label{sphip}
 \end{equation}

\bibliographystyle{siam}
\bibliography{ISGC_refs.bib}
\end{document}